\documentclass[11pt,a4paper]{amsart}
\usepackage[toc,page]{appendix}
\usepackage{graphicx} 
\usepackage{url}
\usepackage{amsmath} 
\usepackage{amssymb}
\usepackage{bm, bbm}
\usepackage{polynom}
\usepackage[outline]{contour}
\usepackage{xcolor}
\usepackage{pgfplots}
\usepackage{mathtools}
\usepackage{stmaryrd}
\usepackage{subcaption}
\pgfplotsset{
  compat=1.13,
}
\usepackage{xfrac}  
\usepackage{enumitem} 
\usepackage{float}
\usepackage{comment}
\usepackage{mathrsfs}
\usepackage[automark]{scrlayer-scrpage}
\usepackage[ruled,vlined, algo2e]{algorithm2e}

\usepackage{tikz}
\usepackage{tikz-cd}
\usepackage{scalerel,stackengine}

\usepackage{natbib}

\theoremstyle{plain}
\newtheorem{thm}{Theorem}[section] 

\newtheorem{rem}[thm]{Remark}

\contourlength{0.1pt}
\contournumber{10}%
\definecolor{clr1}{RGB}{27,158,119}
\definecolor{clr2}{RGB}{217,95,2}
\definecolor{clr3}{RGB}{117,112,179}
\definecolor{clr4}{RGB}{231,41,138}
\definecolor{clr5}{RGB}{102,166,30}
\definecolor{clr6}{RGB}{230,171,2}
\definecolor{clr7}{RGB}{166,118,29}
\definecolor{clr8}{RGB}{102,102,102}

\newcommand{\I}{\mathcal{I}}

\newcommand{\Q}{\mathcal{Q}}

\newcommand{\T}{\mathcal{T}}

\newcommand{\W}{\mathcal{W}}

\newcommand{\calT}{\mathcal{T}}

\newcommand{\LL}{L^2(\Omega)}

\newcommand{\vtheta}{\boldsymbol{\theta}}

\usepackage[most]{tcolorbox}

\definecolor{myGreen2}{RGB}{114,175,30} 
\definecolor{myRed}{RGB}{180,50,50}  
\definecolor{myOrange}{RGB}{225,92,22} 

\definecolor{unia-purple}{RGB}{173, 0, 124}
\definecolor{light-gray}{rgb}{0.8889,0.8889,0.8889}
\definecolor{cpl1}{rgb}{0.8889,0.4356,0.2781}
\definecolor{cpl2}{rgb}{0.0,0.6056,0.9787}
\definecolor{cpl3}{rgb}{0.2422,0.6433,0.3044}

\newlength{\plotwidth}
\newlength{\plotheight}
\usetikzlibrary{arrows.meta}
\usetikzlibrary{backgrounds}
\usepgfplotslibrary{patchplots}
\usepgfplotslibrary{fillbetween}
\pgfplotsset{%
    layers/standard/.define layer set={%
        background,axis background,axis grid,axis ticks,axis lines,axis tick labels,pre main,main,axis descriptions,axis foreground%
    }{
        grid style={/pgfplots/on layer=axis grid},%
        tick style={/pgfplots/on layer=axis ticks},%
        axis line style={/pgfplots/on layer=axis lines},%
        label style={/pgfplots/on layer=axis descriptions},%
        legend style={/pgfplots/on layer=axis descriptions},%
        title style={/pgfplots/on layer=axis descriptions},%
        colorbar style={/pgfplots/on layer=axis descriptions},%
        ticklabel style={/pgfplots/on layer=axis tick labels},%
        axis background@ style={/pgfplots/on layer=axis background},%
        3d box foreground style={/pgfplots/on layer=axis foreground},%
    },
}

\def\dx{\,\text{d}x}

\makeatletter
\newsavebox{\@brx}
\newcommand{\llangle}[1][]{\savebox{\@brx}{\(\m@th{#1\langle}\)}%
  \mathopen{\copy\@brx\mkern2mu\kern-0.9\wd\@brx\usebox{\@brx}}}
\newcommand{\rrangle}[1][]{\savebox{\@brx}{\(\m@th{#1\rangle}\)}%
  \mathclose{\copy\@brx\mkern2mu\kern-0.9\wd\@brx\usebox{\@brx}}}
\makeatother

\newcommand\Nb{\mathtt{N}}

\numberwithin{equation}{section}

\title[Neural Network LOD for Random Coefficients]{Neural Network Localized Orthogonal Decomposition for Numerical Homogenization of Diffusion Operators with Random Coefficients}

\author[F. Kröpfl]{Fabian Kröpfl}
\author[D. Peterseim]{Daniel Peterseim}
\author[E. Ullmann]{Elisabeth Ullmann}
\address[F. Kröpfl]{Institute of Mathematics, University of Augsburg, 
Universit\"atsstra{\ss}e~12a, 86159 Augsburg, Germany}
\email{fabian.kroepfl@uni-a.de}
\address[D. Peterseim]{Institute of Mathematics \& Centre for 
Advanced Analytics and Predictive Sciences (CAAPS), University of 
Augsburg, Universit\"atsstra{\ss}e~12a, 86159 Augsburg, Germany}
\email{daniel.peterseim@uni-a.de}
\address[E. Ullmann]{Department of Mathematics, TUM School of Computation, Information and Technology, Technical University of Munich, Boltzmannstr. 3, 85748 Garching b. München, Germany}
\email{elisabeth.ullmann@tum.de}

\thanks{D.~Peterseim acknowledges support by the Deutsche Forschungsgemeinschaft (DFG, German Research Foundation) through the project~496984632. The work of F.~Kröpfl was part of a project that has received funding from the European Research Council (ERC) under the European Union's Horizon 2020 research and innovation programme (Grant agreement No.~865751 -- RandomMultiScales).}

\begin{document}

\begin{abstract}
This paper presents a neural network--enhanced surrogate modeling approach for diffusion problems with spatially varying random field coefficients. The method builds on numerical homogenization, which compresses fine-scale coefficients into coarse-scale surrogates without requiring periodicity. To overcome computational bottlenecks, we train a neural network to map fine-scale coefficient samples to effective coarse-scale information, enabling the construction of accurate surrogates at the target resolution. This framework allows for the fast and efficient compression of new coefficient realizations, thereby ensuring reliable coarse models and supporting scalable computations for large ensembles of random coefficients.
We demonstrate the efficacy of our approach through systematic numerical experiments for two classes of coefficients, emphasizing the influence of coefficient contrast: (i) lognormal diffusion coefficients, a standard model for uncertain subsurface structures in geophysics, and (ii) hierarchical Gaussian random fields with random correlation lengths.
\end{abstract}

\maketitle

\section{Introduction}
Uncertainty is inherent in the modeling of multiscale diffusion processes in heterogeneous media. In geophysics, materials science, and porous media flow, diffusion coefficients are not only highly oscillatory but also random, reflecting incomplete knowledge of subsurface structures and material properties. Traditional numerical methods struggle to capture the complex interplay between fine-scale heterogeneities and stochastic effects without incurring prohibitive computational costs. Motivated by these challenges, this paper extends the framework introduced in~\cite[Section~6.2]{Kro24} by developing a novel deep-learning-based approach to address the computational bottlenecks in stochastic homogenization.

On a bounded, convex Lipschitz polytope $D\subseteq\mathbb{R}^d$, with $d\in{1,2,3}$, we consider the prototypical diffusion problem%
\begin{equation}\label{eq:modelstrong}
\left\{
\begin{aligned}
-\mathrm{div}(\mathcal{A} \nabla u) &=f \hspace*{0.5cm} \mathrm{in}\ D,\ \\
u &= 0\hspace*{0.5cm} \, \mathrm{on}\ \partial D, \\
\end{aligned} \
\right.
\end{equation}
where the random diffusion coefficient~$\mathcal{A}$ encodes microstructures. We focus on cases where $\mathcal{A}$ is given by a lognormal random field with short correlation lengths. This model is common in geophysics (see, e.g., \cite{Chr13}), where uncertainties in subsurface soil and rock formations arise from prohibitive measurement costs.
Beyond lognormal models, we also consider settings in which the covariance parameters are random, giving rise to hierarchical Gaussian random fields (also called parameterized Gaussian measures) \cite{Latz2019}; see also \cite{Kressner2020}.

The interplay of randomness and multi-scale features renders the approximation of the solution $u$ in~\eqref{eq:modelstrong} particularly challenging. Monte Carlo methods require solving many instances of a deterministic multi-scale PDE with non-periodic diffusion coefficients, and standard mesh-based techniques must resolve all fine-scale oscillations. This full resolution becomes prohibitively expensive when the oscillation length is very small relative to the computational domain.

To mitigate these challenges, numerous numerical homogenization methods have been proposed. Deterministic homogenization has been extensively studied (see, e.g., the recent textbooks and reviews \cite{AHP21, BLeB23, CEH23, MalP20, OS19}), whereas stochastic homogenization has received comparatively less attention, despite several promising approaches. An early overview of computational methods for stochastic homogenization is provided in the review article \cite{ACL11}, which includes an MsFEM-based method for stochastically perturbed regimes, where random fluctuations are modeled as perturbations of an otherwise structured medium.

Representative Volume Element (RVE) methods, which offer a more direct numerical realization of homogenization theory and have a solid theoretical foundation in \cite{Gloria2014, GloriaNeukammOttoPreprint, Gloria2011, Gloria2012}, have garnered increasing attention \cite{BoP04, Cancs2015, GloHa16, Gloria2014, GloNo16, Khoromskaia2020, Mourrat2018}. In this context, variance reduction techniques \cite{Blanc2016, LeBrisLegollMinvielle+2016+25+54} can significantly reduce computational costs; see also the theoretical analysis in \cite{Fischer2019}.

A different approach to numerical stochastic homogenization, which does not rely on classical homogenization theory but instead exploits optimal subspaces with localized basis functions, is based on the Localized Orthogonal Decomposition (LOD) method \cite{GaP19, HeP13, MaP14}. By reformulating LOD via a quasi-local discrete integral operator \cite{GaP17}, a deterministic surrogate model can be constructed that approximates the expected solution on a coarse scale. Despite its favorable theoretical properties, including a rigorous a priori error analysis in the context of quantitative stochastic homogenization theory \cite{FisGP19ppt}, the LOD-based computation of this surrogate model still requires solving fine-scale auxiliary problems on overlapping subdomains for each sample of the random coefficient. This poses a significant computational bottleneck. Methodological improvements on the spatial discretization side, such as the Super-Localized Orthogonal Decomposition (SLOD) \cite{Bonizzoni-Freese-Peterseim, BHP22, pumslod, Freese-Hauck-Peterseim, HaPe21b} and its collocation-type formulation \cite{HMP24}, have enhanced computational efficiency by enabling independent sampling on each patch and better parallelization. However, these improvements do not mitigate the costs associated with sampling.

To address this limitation, we propose a deep-learning-based extension of numerical homogenization. Building on the neural network framework developed in \cite{KroMP21, KroMP23} for compressing deterministic multiscale PDE operators and further elaborated in \cite[Section~6.2]{Kro24}, we extend the approach to the stochastic setting. In this framework, a neural network is trained offline to approximate the mapping from local fine-scale coefficients to effective surrogates. In the subsequent online phase, new coefficient samples can be processed by a single forward pass through the network, drastically reducing computational costs. Moreover, our approach explicitly accounts for high-contrast coefficient realizations, which arise naturally when sampling from lognormal distributions and, despite their practical importance, have received little attention in the existing literature.

This work is structured as follows. In Section~\ref{sec:modelproblem}, we introduce the model problem~\eqref{eq:modelstrong} in its variational formulation and describe the family of random fields considered in this study. Section~\ref{sec:LOD} presents the key principles of the LOD for computing effective numerical surrogate models on a coarse scale for the considered diffusion fields. In Section~\ref{sec:nns}, we summarize the deep-learning-based framework for compressing diffusion fields into numerical surrogates. The proposed method is then investigated through a series of numerical experiments in Section~\ref{sec:numexp:log} and in Section~\ref{sec:numexp:hg}. Finally, Section~\ref{sec:conclusion} concludes the paper and outlines open questions for future research.

\section{Model Problem}\label{sec:modelproblem}
Let $(\Omega,\mathcal{F},\mathbb{P})$ be a suitable probability space, and let the diffusion coefficient $\mathcal{A}: \Omega \times D \to \mathbb{R}$ be a random field. We assume that for $\mathbb{P}$-a.e.\ $\omega\in\Omega$, the coefficient realization $\mathcal{A}(\omega,\cdot)$ is essentially bounded and strictly positive, i.e.,
\begin{equation}\label{eq:pathellipt}
0 < \operatorname*{ess\,inf}_{x\in D} \mathcal{A}(\omega,x) 
\;\leq\; \mathcal{A}(\omega,x) 
\;\leq\; \operatorname*{ess\,sup}_{x\in D} \mathcal{A}(\omega,x) < \infty
\quad \text{for a.e. } x\in D.
\end{equation}
The \emph{contrast} of a realization $\mathcal{A}(\omega)$ is then defined as
\[
\text{contrast}(\omega) \;=\;
\frac{\operatorname*{ess\,sup}_{x\in D}\mathcal{A}(\omega,x)}{\operatorname*{ess\,inf}_{x\in D}\mathcal{A}(\omega,x)}.
\]
Large contrast values correspond to high-contrast coefficients, which are particularly challenging for numerical methods.

Given a deterministic right-hand side $f\in L^2(D)$, we are interested in finding a random field $u: \Omega \times D \rightarrow \mathbb{R}$ that $\mathbb{P}$-a.s. solves problem~\eqref{eq:modelstrong} in a pathwise variational sense. More precisely, we seek a random field $u$ with realizations $u(\omega,\cdot)\in H^1_0(D)$ such that for $\mathbb{P}$-a.e. $\omega\in\Omega$ it holds
\begin{equation}\label{eq:stochweak}
a_{\mathcal{A}(\omega)}(u(\omega),v):= \int_D \mathcal{A}(\omega) \nabla u(\omega) \cdot \nabla v\, \mathrm{d}x = \int_D fv\, \mathrm{d}x  =: \langle f,v \rangle\quad \text{for all } v\in H^1_0(D),
\end{equation}
where the notation $\mathcal{A}(\omega)$ and $u(\omega)$ is motivated by the fact that $\mathcal{A}$ and $u$ can be interpreted as random variables mapping from the space $\Omega$ to $L^\infty(D)$ and $H^1_0(D)$, respectively. 

In this work, we consider lognormal random fields $\mathcal{A}$ of the form $\mathcal{A} = \exp(Z)$, where $Z$ is a centered Gaussian random field with Whittle--Mat\'ern covariance function given by 
\small
\begin{equation}\label{c}
c (x,x') = \frac{\sigma^2}{2^{\nu-1}\Gamma(\nu)} \left(\frac{\sqrt{2\nu}}{\kappa} \|x-x'\|_2\right)^\nu K_\nu\left(\frac{\sqrt{2\nu}}{\kappa} \|x-x'\|_2 \right),\quad \mathrm{for}\, x,\, x'\in D,\ \kappa>0,\ \nu>0.
\end{equation}
\normalsize
In \eqref{c}, $\Gamma(\cdot)$ denotes the Gamma function and $K_\nu(\cdot)$ the modified Bessel function of the second kind and order $\nu$. 
The parameter $\nu$ controls the smoothness of the random field, where smaller values of $\nu$ correspond to rougher realizations, and the parameter $\kappa$ is a characteristic length scale, the aforementioned correlation length. 
The parameter $\sigma^2$ is the variance of the random field $Z$.
Note that for lognormal random fields of the form $\mathcal{A}=\exp(Z)$, where $Z$ is a centered Gaussian random field with Whittle--Mat\'ern covariance function, condition~\eqref{eq:pathellipt} is fulfilled, and thus the path-wise variational formulation~\eqref{eq:stochweak} $\mathbb{P}$-a.s.~has a unique solution. 
However, the contrast in lognormal random fields is not uniformly bounded: it equals the exponential of the difference between the maximum and minimum of the realization of the underlying Gaussian field $Z$ taken over the domain $D$. The contrast is almost surely finite, but it can become extremely large due to the Gaussian field $Z$ taking values in $(-\infty,\infty)$. While a large variance of the Gaussian field increases the spread of typical values and thus makes a high contrast more likely, even modestly rare realizations may exhibit a very large contrast.

So far, we have described models with fixed (deterministic) covariance parameters $(\nu,\kappa,\sigma^2)$. A natural generalization is to introduce uncertainty at the level of these hyperparameters, leading to the notion of  hierarchical Gaussian random fields, sometimes also called  parameterized Gaussian measures~ \cite{Latz2019}; see also \cite{Kressner2020}.
In such a setting, the random coefficient $\mathcal{A}$ is constructed by the composition of two random functions as follows:
First, we select a hyperparameter space $\Theta \subseteq \mathbb{R}^{n}$ 
and a $\Theta$-valued random variable $Y\colon \Omega \rightarrow \Theta$ with a given probability distribution $\mathbb{P}_Y$.
Second, given a sample $\vtheta \in \Theta$ of $Y$ according to the distribution $\mathbb{P}_Y$, we construct a Gaussian random field $Z(\cdot \vert\vtheta) \colon D \times \Omega \rightarrow \mathbb{R} $ with a parameterized mean function $m(x;\vtheta)=\mathbb{E}[Z(x,\cdot\vert \vtheta)]$ and covariance function $c(x,x';\vtheta)=\mathbb{E}[(Z(x,\cdot\vert \vtheta)-m(x;\vtheta))(Z(x',\cdot \vert \vtheta)-m(x';\vtheta))]$, $x, x' \in D$.
Finally, the PDE coefficient function is defined as $\mathcal{A} = \exp(Z)$, where $Z$ denotes the parameterized Gaussian random field $Z(\cdot \vert \vtheta)$.
We are interested in the setting where the parameterized mean function is the zero function, and the parameterized covariance function is the Whittle--Mat\'ern covariance function introduced in \eqref{c}, with one or several random hyperparameters. 
In the literature, various choices of random hyperparameters have been studied; for example, random correlation length $\theta=\kappa$~\cite{Latz2019, Kressner2020}, random variance $\theta=\sigma^2$~\cite{Latz2019}, random smoothness parameter $\theta=\nu$~\cite{Kressner2020}, or combinations of the above, such as $\vtheta=(\nu,\kappa)$~\cite{Kressner2020}.
In our numerical experiments below, we will numerically investigate the case where the correlation length $\kappa$ is a random variable taking values in a bounded interval $[\underline{\theta},\overline{\theta}]$. 

\section{Surrogate Modeling by Localized Orthogonal Decomposition}\label{sec:LOD}
In many multiscale geophysical problems, one is interested in computing quantities of interest, such as the mean function of the solution of~\eqref{eq:stochweak}. When using a Monte Carlo approach, this task requires compressing numerous deterministic multiscale coefficients into appropriate macroscopic surrogates. The Localized Orthogonal Decomposition (LOD) method is specifically designed to effectively approximate these distributed coarse-scale quantities and is, therefore, our framework of choice. Its construction is briefly outlined below.

The main idea behind the LOD is to split the trial space \(H^1_0(D)\) into two components: a finite-dimensional, coefficient-adapted space with superior approximation properties tailored to the problem at hand, and an infinite-dimensional remainder space. To achieve this, let \(\T_H\) be a uniform Cartesian mesh and \(Q^1(\T_H)\) the standard first-order finite element space of piecewise bilinear polynomials. We then define the \(H^1_0\)-conforming space 
\[
V_H := Q^1(\T_H) \cap H^1_0(D).
\]
Based on this setting, we introduce a projective quasi-interpolation operator \(\I_H\colon H^1_0(D) \to V_H\) that satisfies the following stability and approximation properties:
\[
\|H^{-1}(v-\I_H v)\|_{L^2(T)} + \|\nabla \I_H v\|_{L^2(T)} \leq C\, \|\nabla v\|_{L^2(\Nb(T))}
\]
for all \(v\in H^1_0(D)\) and any element \(T\in \T_H\). Here, the constant \(C\) is independent of \(H\), and \(\Nb(S) := \Nb^1(S)\) denotes the \emph{element neighborhood} of a subset \(S\subseteq D\), defined by
\[
\Nb^1(S) := \bigcup \Bigl\{\overline{T} \in \T_H \,\Big|\, \overline{S}\cap \overline{T}\neq \emptyset\Bigr\}.
\]
The operator $\I_H$ now allows us to characterize a so-called \emph{fine-scale space} $\W$ as the space of all functions in $H^1_0(D)$ which cannot be captured well within the finite element space $V_H$, i.e., we set
\begin{equation*}
\W := \ker \I_H\vert_{H^1_0(D)}.
\end{equation*} 
Moreover, for any set $S \subseteq D$, we define a local fine-scale space
\begin{equation*}
\W(S) := \{w \in \W\;\vert\;\mathrm{supp}(w) \subseteq S\}.
\end{equation*}
For a given realization \(A := \mathcal{A}(\omega)\) of the random field \(\mathcal{A}\), the goal is to construct a multi-scale trial space by \emph{correcting} the nodal basis functions of \(V_H\) with suitably chosen fine-scale functions from \(\W\). The resulting corrected basis functions span a new multi-scale space of the same dimension as \(V_H\), but with significantly improved approximation properties due to the fine-scale information encoded in the basis functions.

To this end, for \(\ell \geq 2\) we first define the \emph{element neighborhood of order \(\ell\)} as
\[
\Nb^\ell(S) := \Nb\bigl(\Nb^{\ell-1}(S)\bigr),
\]
and for a given function \(v_H \in V_H\), the \emph{element corrector} \(\Q_{A,T}^\ell v_H \in \W(\Nb^\ell(T))\) corresponding to the element \(T \in \T_H\) is the unique solution of
\begin{equation}\label{eq:corT}
a_{A}\bigl(\Q^\ell_{A,T} v_H, w\bigr) = \int_T A \nabla v_H \cdot \nabla w \, dx \quad \text{for all } w \in \W(\Nb^\ell(T)).
\end{equation}
Note that the support of the element corrector $Q_{A,T}^\ell v_H$ is limited to the element neighborhood $\Nb^\ell(T)$ by construction, and in practice, problem~\eqref{eq:corT} has to be approximated on a finer discretization scale that sufficiently resolves all fine-scale oscillations of the realization $A$. The \emph{global corrector} $Q_{A}^\ell v_H\in\W$ is then obtained by summing the individual element correctors, i.e., by computing
\begin{equation}\label{eq: corglob}
Q_{A}^\ell v_H := \sum_{T\in\T_H} \Q_{A,T}^\ell v_H.
\end{equation}
The computations of the individual element correctors are mutually  independent and only depend on the values of the realization $A$ in the respective element neighborhoods. Global computations on a fine discretization scale can thus be avoided.
If one chooses $\ell = \infty$, however, which corresponds to a computation of the element correctors over the entire computational domain, one obtains the orthogonality relation
\begin{equation}\label{eq: orthdecomp}
a_{A}((\mathsf{id}-\Q^\infty_{A}) v_H, w) = 0\quad \text{for all } w \in \W,
\end{equation}
where $\Q^\infty_{A}: V_H \to \W$ denotes the correction operator mapping a function in $V_H$ to its global corrector. Note that the relation~\eqref{eq: orthdecomp} induces the $a_A$-orthogonal splitting $H^1_0(D) = (\mathsf{id} - \Q_A^\infty)V_H \oplus \W$, with an ``ideal'' multi-scale trial space $(\mathsf{id} - \Q_A^\infty)V_H$. For the aforementioned reasons, it is, however, very impractical to use this space in the actual implementation of the method. However, it can be shown that the element correctors $\Q_{A,T}^\infty v_H$ that contribute to the global corrector decay exponentially fast away from $T$, which motivates the use of only localized correctors and the corresponding localized multi-scale space $(\mathsf{id} - \Q_A^\ell)V_H$ in practice. 

For more details on the LOD methodology and a proof of the localization result, we refer the reader to the original work \cite{MaP14} and the textbook \cite{MalP20}. Alternative derivations and generalizations in the spirit of gamblets \cite{OS19} are discussed in the review \cite{AHP21}. However, these approaches are often not easily adapted to the stochastic setting, as they typically do not support a Petrov--Galerkin formulation that keeps the test functions deterministic in a straightforward manner. This feature is crucial in our context to avoid the need for sampling inner products of two random basis functions. The same limitation applies to Galerkin versions of super-localized methods, with the exception of the collocation variant \cite{HMP24}.

Using this construction, the Petrov--Galerkin LOD (PG-LOD) for approximating~\eqref{eq:stochweak} in a macroscopic sense now seeks $u_H\in V_H$ such that
\begin{equation}\label{eq:discvarproblempg}
a_A((\mathsf{id}-\Q^\ell_A) u_H, v_H) = \langle f,\, v_H \rangle \quad\text{for all } v_H \in V_H,
\end{equation}
which is a first-order accurate approximation in $L^2(D)$ if one chooses $\ell \approx \vert\log H\vert$. Given an enumeration of the degrees of freedom $\mathcal{N}(\T_H)$ and the corresponding nodal basis functions $\Lambda_j$ of $V_H$, the global PG-LOD system matrix $\mathbf{S}_A$ is given by
\begin{equation}\label{eq:Sglobal}
(\mathbf{S}_A)_{i,j} = a_A((\mathsf{id}-\Q^\ell_A) \Lambda_j, \Lambda_i),\quad i,j \in \mathcal{N}(\T_H),
\end{equation}
then serves as the desired macroscopic surrogate for the realization $A$.

In practice, the global matrix $\mathbf{S}_A$ is assembled from local contributions $\mathbf{S}_{A,T}$ corresponding to the element neighborhoods around the macroscopic elements $T\in\T_H$. For a given element $T\in\T_H$, this matrix is defined by
\begin{equation}\label{eq:Slocal}
(\mathbf{S}_{A,T})_{i,j} = \int_{\Nb^\ell(T)}A\, \nabla(\mathsf{id}-\Q^\ell_A)\Lambda_j\vert_T \cdot \nabla\Lambda_i\, \dx,\quad j\in\mathcal{N}(T),\ i\in \mathcal{N}(\Nb^\ell(T)). 
\end{equation}
Similar to classical finite element assembly processes, assembling $\mathbf{S}_A$ from the local contributions $\mathbf{S}_{A,T}$ can be realized by using appropriate local-to-global mappings $\Phi_T$ that transform the local matrix for given $T \in \calT_H$ into an inflated matrix with respect to the global degrees of freedom. The global stiffness matrix $\mathbf{S}_A$ can thus be decomposed as
\begin{equation}\label{eq:decSLOD}
\mathbf{S}_A = \sum_{T \in \T_H} \Phi_T\big(\mathbf{S}_{A,T}\big),
\end{equation}
see~\cite[Sec.~3.4]{KroMP21} for more details. As mentioned above, the local matrices $\mathbf{S}_{A,T}$ depend on $A$, but the local-to-global mappings $\Phi_T$ are independent of any given realization.
Once the global matrix has been assembled, it can then be used to reliably approximate the solution to~\eqref{eq:stochweak} on the chosen target scale.

\section{Neural Network Compression}\label{sec:nns}
This section briefly summarizes the main ideas from the framework proposed in~\cite{KroMP21}. For more details, we refer to this article or~\cite{Kro24}. 
The main idea of this framework is to approximate the local contributions to the PG-LOD stiffness matrix, i.e., the matrices $\mathbf{S}_{A,T}$ introduced in~\eqref{eq:Slocal}, by a single feedforward neural network of moderate size.

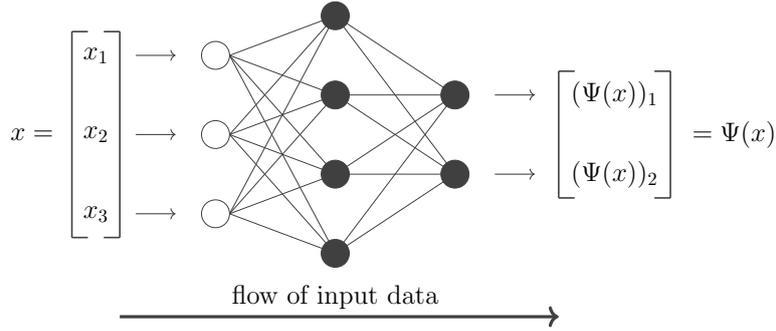
\begin{figure}
\begin{center}
\scalebox{0.7}{\begin{tikzpicture}[scale=1.5]

\node at (-2.3,2) {\Large $x = $};
\node at (-1.5,1) {\Large $x_3$};
\node at (-1.5,2) {\Large $x_2$};
\node at (-1.5,3) {\Large $x_1$};

\draw[darkgray,line width=0.3mm] (-1.8,0.7) -- (-1.8,3.3);
\draw[darkgray,line width=0.3mm] (-1.8,0.7) -- (-1.6,0.7);
\draw[darkgray,line width=0.3mm] (-1.8,3.3) -- (-1.6,3.3);

\draw[darkgray,line width=0.3mm] (-1.2,0.7) -- (-1.2,3.3);
\draw[darkgray,line width=0.3mm] (-1.2,0.7) -- (-1.4,0.7);
\draw[darkgray,line width=0.3mm] (-1.2,3.3) -- (-1.4,3.3);

\filldraw[->,line width=0.2mm,darkgray] (-1,1) -- (-0.5,1);
\filldraw[->,line width=0.2mm,darkgray] (-1,2) -- (-0.5,2);
\filldraw[->,line width=0.2mm,darkgray] (-1,3) -- (-0.5,3);

\draw[darkgray]  (0,1) circle (5pt);
\draw[darkgray]  (0,2) circle (5pt);
\draw[darkgray]  (0,3) circle (5pt); 

\filldraw[darkgray]  (1.5,0.5) circle (5pt);
\filldraw[darkgray]  (1.5,1.5) circle (5pt);
\filldraw[darkgray]  (1.5,2.5) circle (5pt);
\filldraw[darkgray]  (1.5,3.5) circle (5pt);

\filldraw[darkgray]  (3,1.5) circle (5pt);
\filldraw[darkgray]  (3,2.5) circle (5pt);

\draw[darkgray] (1 + 4pt,1) -- (1.5 ,0.5);
\draw[darkgray] (1 + 4pt,1) -- (1.5 ,1.5);	
\draw[darkgray] (1 + 4pt,1) -- (1.5 ,2.5);	
\draw[darkgray] (1 + 4pt,1) -- (1.5 ,3.5);	

\draw[darkgray] (1 + 4pt,2) -- (1.5 ,0.5);
\draw[darkgray] (1 + 4pt,2) -- (1.5 ,1.5);	
\draw[darkgray] (1 + 4pt,2) -- (1.5 ,2.5);	
\draw[darkgray] (1 + 4pt,2) -- (1.5 ,3.5);

\draw[darkgray] (1 + 4pt,3) -- (1.5 ,0.5);
\draw[darkgray] (1 + 4pt,3) -- (1.5 ,1.5);	
\draw[darkgray] (1 + 4pt,3) -- (1.5 ,2.5);	
\draw[darkgray] (1 + 4pt,3) -- (1.5 ,3.5);

\draw[darkgray] (1.5 ,0.5) -- (3,1.5);
\draw[darkgray] (1.5 ,0.5) -- (3,2.5);	
	
\draw[darkgray] (1.5 ,1.5) -- (3,1.5);
\draw[darkgray] (1.5 ,1.5) -- (3,2.5);	

\draw[darkgray] (1.5 ,2.5) -- (3,1.5);
\draw[darkgray] (1.5 ,2.5) -- (3,2.5);

\draw[darkgray] (1.5 ,3.5) -- (3,1.5);
\draw[darkgray] (1.5 ,3.5) -- (3,2.5);	

\filldraw[->,line width=0.2mm,darkgray] (3.5,1.5) -- (4,1.5);
\filldraw[->,line width=0.2mm,darkgray] (3.5,2.5) -- (4,2.5);

\node at (5,1.5) {\Large $(\Psi(x))_2$};
\node at (5,2.5) {\Large $(\Psi(x))_1$};
\node at (6.5,2) {\Large $= \Psi(x)$};

\draw[darkgray,line width=0.3mm] (4.3,1.2) -- (4.3,2.8);
\draw[darkgray,line width=0.3mm] (4.3,2.8) -- (4.5,2.8);
\draw[darkgray,line width=0.3mm] (4.3,1.2) -- (4.5,1.2);

\draw[darkgray,line width=0.3mm] (5.7,1.2) -- (5.7,2.8);
\draw[darkgray,line width=0.3mm] (5.5,1.2) -- (5.7,1.2);
\draw[darkgray,line width=0.3mm] (5.5,2.8) -- (5.7,2.8);

\filldraw[->,line width=0.6mm,darkgray] (-1.2,-0.3) -- (4.3,-0.3);
\node at (1.5,-0.05) {\Large flow of input data};

\end{tikzpicture} }
\end{center}
\caption[Illustration of a two-layer neural network.]{Illustration of a feedforward neural network with two layers.}\label{fig:dffnn}
\end{figure}

Feedforward neural networks are nested functions of the form
\begin{equation*}
\Psi(x) = \Psi_L(\Psi_{L-1}(\dots(\Psi_1(x))\dots )),    
\end{equation*}
where $L\in\mathbb{N}$ is the depth and the individual functions 
\begin{equation*}
\Phi_l(x) = \rho_l(W_l x + b_l)    
\end{equation*}
are the layers of the network with the weight matrix $W_l$ and the bias vector $b_l$.
The functions $\rho_l: \mathbb{R}\to\mathbb{R}$ are nonlinear activation functions that act component-wise on vectors by convention. In this work, we will only consider the classic ReLU activation function given by $\rho_l(x) = \max\{0,x\}$ in the hidden layers and the identity function in the last layer. The nested structure of neural networks can be visualized as a computational graph, where the sparsity patterns of the weight matrices determine the connectivity of the individual layers. Figure~\ref{fig:dffnn} (adapted from~\cite{Kro24}) illustrates a fully connected two-layer neural network $\Psi$, i.e., a network with fully occupied weight matrices that takes vectors in $\mathbb{R}^3$ as input and outputs vectors in $\mathbb{R}^2$. 

Since neural networks only take discrete values as input, we first have to discretize the realizations of the considered random field.  Let $R_T: \Omega \to \mathbb{R}^m$ be the operator that maps any $\omega\in\Omega$ to $m$ discrete relevant features of $\mathcal{A}(\omega)$, which then serve as input to the neural network. In our present setting, a canonical choice for the operator $R_T$ is to map any sample $\omega$ to a piece-wise constant approximation of the realization $\mathcal{A}(\omega)$ on a finer local mesh $\T_\varepsilon(\Nb^\ell(T))$ on $\Nb^\ell(T)$ that sufficiently resolves the fine-scale oscillations of $\mathcal{A}(\omega)$. 

Since neural networks have a fixed input dimension, we also have to ensure that for each $\omega \in \Omega$ and each $T \in \T_H$, all discrete input vectors $R_T(\omega)$ are of uniform dimension, regardless of the location of the central element $T$ with respect to the boundary. To this end, we artificially extend the domain $D$ and the mesh $\T_H$ by $\ell$ layers of outer elements around the boundary of $D$. On these outer elements, we assign the value $\mathcal{A}=0$. This padding convention is not a modeling assumption, and the random field $\mathcal{A}$ could also be sampled outside of $D$. However, it provides an easy and consistent way to ensure that every neighborhood $\Nb^\ell(T)$ contains the same number of coarse elements. Assigning zero values outside of $D$ not only ensures a uniform input dimension but also implicitly encodes the global boundary into the local patches, so that the network can distinguish between interior and boundary neighborhoods. For consistency, the finer local mesh $\T_\varepsilon(\Nb^\ell(T))$ is chosen as a uniform refinement of $\T_H(\Nb^\ell(T))$, which ensures that every patch is represented by the same number of fine elements. In this way, we obtain a uniform input dimension of $m = \lvert \T_\varepsilon(\Nb^\ell(T)) \rvert$.
This extension introduces additional nodes outside of the physical domain. Consequently, all local system matrices $\mathbf{S}_{A,T}$ (and hence all network outputs) are of uniform size $N^\ell_H \times 2^d$, where $N^\ell_H = (2\ell + 2)^d$ denotes the number of nodes in $\Nb^\ell(T)$. The rows associated with these outer nodes are disregarded during the assembly of the global system matrix, so the Dirichlet boundary condition is still enforced exactly at the global level.

Let now $\Psi_\mathbf{w}: \mathbb{R}^m\to \mathbb{R}^{N^\ell_H \cdot 2^d}$ be a neural network with parameters $\mathbf{w}=\{W_l,b_l\}_{l=1}^L$ encoding the weights and biases at all layers that map locally discretized realizations $A = \mathcal{A}(\omega)$ of the random field to an approximation of the (column-wise vectorized) local PG-LOD matrix $\mathbf{S}_{A,T}$. 
The method for approximating the global matrix $\mathbf{S}_A$ now retains the decomposition~\eqref{eq:decSLOD} of the global matrix and replaces the local contributions $\mathbf{S}_{A,T}$ with suitable approximations based on the neural network $\Psi$. That is, for a given $\omega\in\Omega$ and the corresponding realization $A = \mathcal{A}(\omega)$, we approximate the global PG-LOD matrix $\mathbf{S}_A$ by its neural network-based surrogate
\begin{equation*}
\widehat{\mathbf{S}}_A = \sum_{T \in \T_H} \Phi_T\big(\Psi_\mathbf{w}(R_T(\omega)\big).
\end{equation*}
This is a reasonable approach because the local corrector problems~\eqref{eq:corT} are dependent on the actual values of the realization, independent of their position within the domain. 

The network is first trained on a dataset consisting of local coefficient-surrogate pairs $(R_T(\omega^{(i)}), \mathbf{S}_{A,T}^{(i)})$ by minimizing the loss function
\begin{equation}\label{eq:loss}
\mathcal{J}(\mathbf{w}) = \frac{1}{2n\cdot \mathrm{card}(\calT_H)} \sum_{i=1}^n \sum_{T\in \calT_H} \ \frac{\| \Psi_\mathbf{w}(R_T(\omega^{(i)})) - \mathrm{vec}(\mathbf{S}_{A,T}^{(i)}) \|^2_2}{\| \mathrm{vec}(\mathbf{S}_{A,T}^{(i)}) \|^2_2},
\end{equation}
over a space of admissible network parameters, i.e., all entries of the network's weight matrices and bias vectors. Here, the $\mathrm{vec}$-operator flattens the system matrices column-wise into vectors.

Our approach has the advantage that, after the initial training phase, the process of compressing new realizations into macroscopic surrogates is substantially accelerated compared to the usual PG-LOD algorithm. This is due to the fact that, instead of computing the local contributions by solving corrector problems of the form~\eqref{eq:corT}, only a forward pass of the mini-batch $(R_T(\omega))_{T\in\T_H}$ through the network needs to be performed.

\section{Numerical Experiments with lognormal diffusion coefficients}\label{sec:numexp:log}

In this section, we investigate the feasibility of the neural-network-based compression approach described in Section~\ref{sec:nns} for lognormal random fields, as outlined in Section~\ref{sec:modelproblem}.
Specifically, we explore how the contrast of the random diffusion coefficient affects the training of a given network architecture in terms of its ability to fit a given dataset. 
From a theoretical point of view, increasing the contrast of the diffusion coefficient leads to increased bounds on the depth and number of non-zero parameters of an approximating neural network; see~\cite{KroMP23}. For a fixed network depth and number of network parameters, as in our test problem, we therefore expect that the best possible training and validation errors that can be achieved are small for low contrast coefficients and increase as the contrast of the coefficients increases. 

\subsection{Problem parameters and neural network}\label{sec:setup}

In all experiments, we consider the unit square $D=(0,1)^2$ as our computational domain. 
The discretization resp. problem parameters are $H=2^{-4}$, $\varepsilon = 2^{-6}$, $h = 2^{-9}$, $\ell = 2$, and $f\equiv 1$.
We use the same network as in the numerical experiments in~\cite{KroMP21}, i.e., a feedforward neural network with eight layers and ReLU activation function. The weight matrices and bias vectors have the dimensions
\begin{equation*}
\begin{aligned}
&\mathbf{W}_{1} \in \mathbb{R}^{1600\times 1600}, \hspace*{0.1cm} &&\mathbf{W}_{2} \in \mathbb{R}^{800\times 1600}, \hspace*{0.1cm}
&&\mathbf{W}_{3} \in \mathbb{R}^{800\times 800}, \hspace*{0.1cm} &&\mathbf{W}_{4} \in \mathbb{R}^{400\times 800}, \\ 
&\mathbf{W}_{5} \in \mathbb{R}^{400\times 400}, \hspace*{0.1cm} &&\mathbf{W}_{6} \in \mathbb{R}^{144\times 400}, \hspace*{0.1cm}
&&\mathbf{W}_{7} \in \mathbb{R}^{144\times 144}, \hspace*{0.1cm} &&\mathbf{W}_{8} \in \mathbb{R}^{144\times 144}, \\
&\, \mathbf{b}_{1} \in \mathbb{R}^{1600}, \hspace*{0.1cm} &&\,\mathbf{b}_{2}, \in \mathbb{R}^{800}, \hspace*{0.1cm}
&&\,\mathbf{b}_{3} \in \mathbb{R}^{800}, \hspace*{0.1cm} &&\,\mathbf{b}_{4} \in \mathbb{R}^{400}, \\ 
&\,\mathbf{b}_{5} \in \mathbb{R}^{400}, \hspace*{0.1cm} &&\,\mathbf{b}_{6} \in \mathbb{R}^{144}, \hspace*{0.1cm}
&&\,\mathbf{b}_{7} \in \mathbb{R}^{144}, \hspace*{0.1cm} &&\,\mathbf{b}_{8} \in \mathbb{R}^{144}.
\end{aligned}
\end{equation*}
Instead of training the network from scratch, starting from a randomly initialized set of weight matrices and bias vectors, we employ transfer learning: We reuse the weights and biases obtained in the numerical experiments in~\cite{KroMP21}, where this architecture was trained on a family of piecewise constant, uniformly elliptic coefficients.
We denote this network as $\widehat{\Psi}^{\mathrm{pg}}$.

\subsection{Lognormal diffusion coefficient}\label{sec:lognormal}

The diffusion coefficient $\mathcal{A}=\exp(Z)$ in \eqref{eq:modelstrong} is a lognormal random field.
We consider different variance parameters $\sigma^2$ in the covariance function of $Z$ in \eqref{c}. 
Let $c_{\sigma^2}$ denote the Whittle--Matérn covariance function with variance $\sigma^2$. Since we are mainly interested in exploring the effect of increasing the contrast on the deep learning algorithm, and increasing the contrast of a typical realization of a lognormal random field amounts to increasing the variance of the underlying Gaussian random field, we keep the other parameters of $c_{\sigma^2}$ fixed for now. The smoothness parameter is fixed at the value $\nu:=1$ throughout all experiments below, and the (for now) deterministic correlation length is $\kappa := 2^{-6}$. Accordingly, we denote a centered Gaussian random field with covariance function $c_{\sigma^2}$ by $Z_{\sigma^2}.$
We interpret a lognormal random field as a family of deterministic coefficients that is parameterized by $\omega \in \Omega$; that is, we consider coefficient classes of the form
\begin{equation}
\mathcal{A}_{\sigma^2} := \{ \exp(Z_{\sigma^2}(\omega))\ \vert\ \omega\in\Omega \},
\end{equation}
which are indexed by the variance parameter $\sigma^2$.
In the following, we consider three situations: namely, the low contrast case with $\sigma^2=0.5$, the moderate-contrast case with $\sigma^2=1$, and the high-contrast case with $\sigma^2=2$.
Note that the contrasts of all three random fields are not uniformly bounded in $\omega$. However, since the density function of a centered Gaussian random variable decays exponentially fast away from zero, the probability of obtaining high-contrast realizations when drawing samples from a lognormal random field is close to zero for small variance parameters $\sigma^2$ and increases steadily as $\sigma^2$ is increased. Since typical realizations of $\mathcal{A}_{0.5}$, $\mathcal{A}_{1}$, and $\mathcal{A}_2$ fall into the comparatively low-contrast, moderate-contrast, and high-contrast regimes, respectively, we will sometimes refer to $\mathcal{A}_{0.5}$ as the \emph{low-contrast random field}, $\mathcal{A}_{1}$ \emph{as the moderate-contrast random field} , and $\mathcal{A}_{2}$ as the \emph{high-contrast random field}. Analogously, we will use the terms \emph{low-contrast realizations, moderate-contrast realizations} , and \emph{high-contrast realizations} for the realizations of those random fields, respectively.

We generate six datasets with increasing variance parameters: two for the low-contrast case $\sigma^2 = 0.5$, two for the moderate-contrast case $\sigma^2 = 1$, and two for the high-contrast case $\sigma^2 = 2$, corresponding to the coefficient classes $\mathcal{A}_{0.5}$, $\mathcal{A}_{1}$ and $\mathcal{A}_2$, respectively. In the following, we denote these datasets by $\mathfrak{d}_{0.5},\mathfrak{D}_{0.5}, \mathfrak{d}_1, \mathfrak{D}_1, \mathfrak{d}_2$ and  $\mathfrak{D}_2$, where $\mathfrak{d}$ refers to a small dataset and $\mathfrak{D}$ refers to a large dataset. 

The number of random field realizations in each dataset is given in Table~\ref{tab:size}.
\renewcommand{\arraystretch}{1.2}
\begin{table}[ht]
\begin{tabular}{l|cccccc}\hline
Complete dataset& $\mathfrak{d}_{0.5}$ & $\mathfrak{d}_{1}$ & $\mathfrak{d}_{2}$ & $\mathfrak{D}_{0.5}$ & $\mathfrak{D}_{1}$ & $\mathfrak{D}_{2}$ \\
Number of realizations & 300 & 300 & 300 & 800 & 3,200 & 3,200\\\hline
\end{tabular}
\caption{The size of the datasets.}\label{tab:size}
\end{table}
In order to generate training, validation, and test data for our experiments, for every $\sigma^2\in\{0.5,1,2\}$, we first generate $300$ (resp. $800$ or $3200$ for the large datasets) global realizations $Z_{\sigma^2}^{(i)}, i = 1,\dots,300/800/3200$ of the Gaussian random field $Z_{\sigma^2}$ at the midpoints of a uniform Cartesian mesh $\calT_\varepsilon$ of $D=[0,1]^2$ with a mesh size $\varepsilon = 2^{-7}$. For the generation of the realizations, we use the so-called \emph{Circulant Embedding} method introduced in~\cite{DieN97}; see also \cite[Ch.~7]{LPS2014}.
We use the MATLAB code for the Circulant Embedding method provided in~\cite{KroB14}, which we have adapted to meet our needs.

\subsection{Neural network inputs and outputs}

We split each sample $Z^{(i)}_{\sigma^2}$ of the Gaussian random field into sub-samples $Z_{\sigma^2,T}^{(i)}$ by applying the operators $R_T$. The splitting is based on {element neighborhoods} $\Nb^\ell(T)$ for $\ell = 2$, which are centered around the macroscopic elements $T\in\calT_H$. Since we work on a uniform Cartesian mesh at the target scale $H = 2^{-4}$, we obtain $256$ sub-samples $Z_{\sigma^2,T}^{(i)} \in\mathbb{R}^{1600}$ per global realization $Z_{\sigma^2}^{(i)}\in\mathbb{R}^{16384}$.
Then we compute the exponential of all these sub-samples, i.e., we compute $\mathcal{A}_{\sigma^2,T}^{(i)}:= \exp(Z_{\sigma^2,T}^{(i)})$ for every $\sigma^2$, $i$, and $T$, in the sense that the exponential function is applied component-wise to every entry of $Z_{\sigma^2,T}^{(i)}\in\mathbb{R}^{1600}.$ 
By construction, all local sub-samples $\mathcal{A}_{\sigma^2,T}^{(i)}$ take the value zero on the artificially introduced elements outside of the physical domain $D$.
Hence, we set the entries corresponding to those elements equal to zero.
We then compute the local surrogate matrices $(\mathbf{S}_{\mathcal{A}_{\sigma^2},T}^{\mathrm{pg}})^{(i)}\in\mathbb{R}^{36\times 4}$ with the PG-LOD, where the corrector problems are approximated on a fine mesh of mesh size $h = 2^{-9}$.
Then, we flatten the surrogate matrices column-wise into vectors in $\mathbb{R}^{144}$.

We assign the realizations in Table~\ref{tab:size} into training, validation, and test sets $\mathfrak{d}_{\sigma^2}^{\mathrm{train}},\mathfrak{d}_{\sigma^2}^{\mathrm{val}}$, and $\mathfrak{d}_{\sigma^2}^{\mathrm{test}}$ (resp. $\mathfrak{D}_{\sigma^2}^{\mathrm{train}},\mathfrak{D}_{\sigma^2}^{\mathrm{val}}$ and $\mathfrak{D}_{\sigma^2}^{\mathrm{test}}$) according to an $80:10:10$ split. 
That is, the first $80\%$ of the realizations are assigned to the training set, the next $10 \%$ are assigned to the validation set, and the last $10\%$ are assigned to the test set.
This gives the following numbers:
For every $\sigma\in\{0.5,1,2\}$, we obtain $240\,\cdot\, 256 = 61,440$ input-output pairs $ (Z_{{\sigma^2},T}^{(i)} ,(\mathbf{S}_{\mathcal{A}_{\sigma^2},T}^{\mathrm{pg}})^{(i)})$ for the training set $\mathfrak{d}_{\sigma^2}^{\mathrm{train}}$ and $30 \,\cdot\, 256 = 7,680$ pairs for the validation set $\mathfrak{d}_{\sigma^2}^{\mathrm{val}}$ and the test set $\mathfrak{d}_{\sigma^2}^{\mathrm{test}}$, respectively. 
For the large, low-contrast dataset we obtain
$640 \,\cdot\, 256$ pairs in $\mathfrak{D}_{0.5}^{\mathrm{train}}$ and $80\,\cdot\, 256$ pairs in $\mathfrak{D}_{0.5}^{\mathrm{val}}$ and  $\mathfrak{D}_{0.5}^{\mathrm{test}}$. For the large, medium- and high-contrast datasets we obtain $2560\,\cdot\, 256$ pairs in $\mathfrak{D}_{1}^{\mathrm{train}}$ and $\mathfrak{D}_{2}^{\mathrm{train}}$, and $320\,\cdot\, 256$ pairs in the associated validation and test sets.
It is important to note that we use the local sub-samples $Z_{{\sigma^2},T}^{(i)}$ corresponding to the realizations of the Gaussian random field $Z_{\sigma^2}$ as input to our neural network. 
Hence, our neural network approximates the mapping $Z_{{\sigma^2},T} \mapsto \mathbf{S}_{\mathcal{A}_{\sigma^2},T}^{\mathrm{pg}}$.
\begin{rem}
In experiments not reported here, we trained the neural network based on the inputs $\mathcal{A}_{{\sigma^2},T}^{(i)}$.
That is, we approximated the mapping $\mathcal{A}_{{\sigma^2},T} \mapsto \mathbf{S}_{\mathcal{A}_{\sigma^2},T}^{\mathrm{pg}}$.
However, we abandoned this approach since it was much more straightforward to obtain good results using the Gaussian inputs. 
A possible explanation is that the Gaussian inputs are centered around zero and show much less variation in terms of their order of magnitude.
\end{rem}

\subsection{Neural network training}
The starting point is the neural network $\widehat{\Psi}^{\mathrm{pg}}$ introduced in section~\ref{sec:setup}.
In our numerical experiments in sections~\ref{subsec:loco}--\ref{subsec:hico}, we retrain the network $\widehat{\Psi}^{\mathrm{pg}}$ separately on the datasets $\mathfrak{d}_{0.5}$ vs. $\mathfrak{D}_{0.5}$, $\mathfrak{d}_1$ vs. $\mathfrak{D}_1$, and $\mathfrak{d}_2$ vs. $\mathfrak{D}_2$ to study the impact of increasing contrast and the size of the dataset in a systematic way. 
For the small datasets, in each of the three experiments below, we train the network $\widehat{\Psi}^{\mathrm{pg}}$ for $60$ epochs on the respective training sets $\mathfrak{d}_{\sigma^2}^{\mathrm{train}}$, while also monitoring the loss on the corresponding validation sets $\mathfrak{d}_{\sigma^2}^{\mathrm{val}}$. 
For the large datasets, the network $\widehat{\Psi}^{\mathrm{pg}}$ is retrained for $40$ epochs on the respective training sets $\mathfrak{D}_{\sigma^2}^{\mathrm{train}}$, while monitoring $\mathfrak{D}_{\sigma^2}^{\mathrm{val}}$. 
The networks obtained after retraining are denoted by $\widehat{\Psi}_{\sigma^2,\mathfrak{d}}^{\mathrm{pg}}$ and $\widehat{\Psi}_{\sigma^2,\mathfrak{D}}^{\mathrm{pg}}$, respectively. 
In all cases, we use a batch size of $100$ and the Adam optimizer with a step size of $10^{-3}$ for the first $30$ epochs, a step size of $10^{-4}$ for the epochs $31$ to $40$, and a step size of $9.5\cdot 10^{-5}$ for the epochs $41$ to $60$. We experimented extensively with different configurations of step size strategies, optimization algorithms, and batch sizes; this particular configuration yielded the best results across all considered datasets.
After retraining the network $\widehat{\Psi}^{\mathrm{pg}}$ on the respective dataset, we always generate one additional realization of $\mathcal{A}_{\sigma^2}$, hereafter denoted by $\mathcal{A}_{\sigma^2}^{(t)}$. For this additional realization, we compute the corresponding coarse-scale surrogate $(\mathbf{S}_{\mathcal{A}_{\sigma^2}}^{\mathrm{pg}})^{(t)}$ based on the PG-LOD and its neural-network-based approximation $(\widehat{\mathbf{S}}_{\mathcal{A}_{\sigma^2}}^{\mathrm{nn}})^{(t)}$ based on the retrained network $\widehat{\Psi}^{\mathrm{pg}}_{\sigma^2}$.
We then compare the approximate solutions of the model problem~\eqref{eq:stochweak} for the diffusion coefficient $\mathcal{A}_{\sigma^2}^{(t)}$ and the right-hand side $f\equiv 1$.
We describe the results below.\\

\subsection{Experiment 1: Low-contrast realizations}~\label{subsec:loco}
\begin{figure}
     \centering
     \begin{subfigure}[a]{0.49\textwidth}
         \centering
         \includegraphics[width=\textwidth,trim= 80 0 100 0,clip]{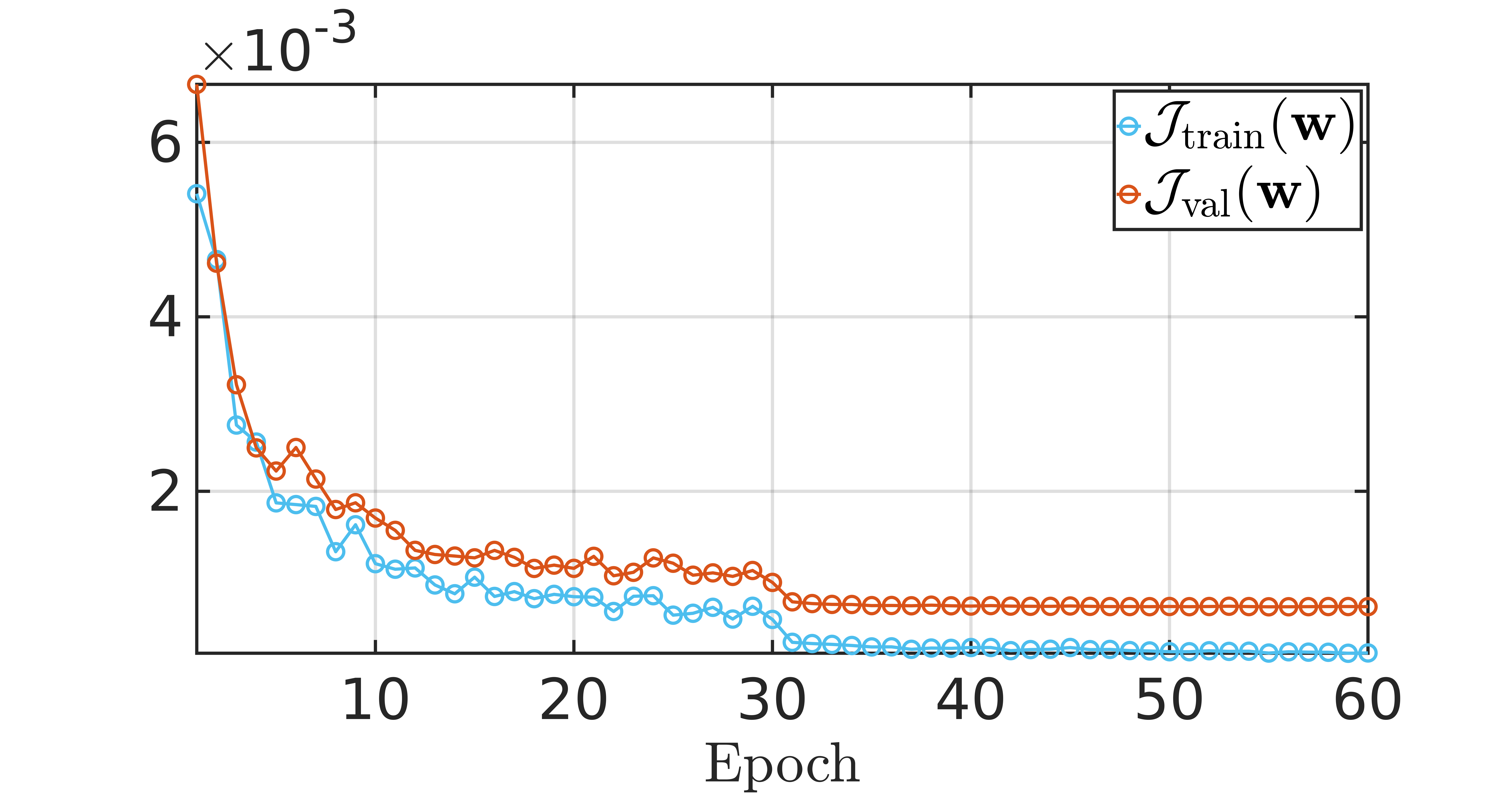}
     \end{subfigure} 
          \begin{subfigure}[a]{0.49\textwidth}
         \centering
         \includegraphics[width=\textwidth,trim= 80 0 100 0,clip]{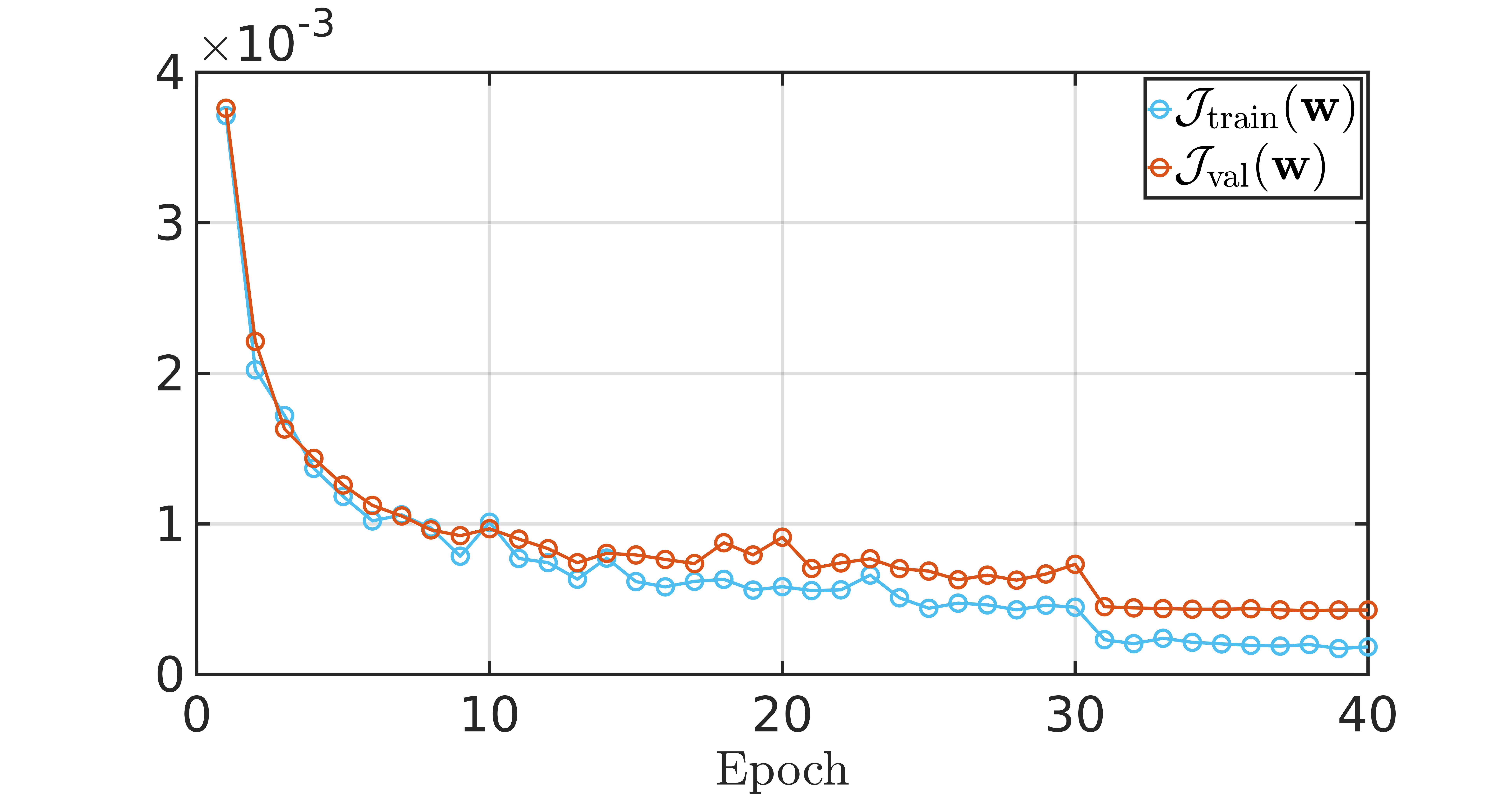}
     \end{subfigure} 
\caption[Development of loss for low-contrast dataset.]{Experiment 1: Training and validation loss for the low-contrast datasets~$\mathfrak{d}_{0.5}$ (left) and~$\mathfrak{D}_{0.5}$ (right).}
\label{fig:lossloco}
\end{figure}

\begin{figure}
     \centering
     \begin{subfigure}[a]{0.45\textwidth}
         \centering
         \includegraphics[width=\textwidth]{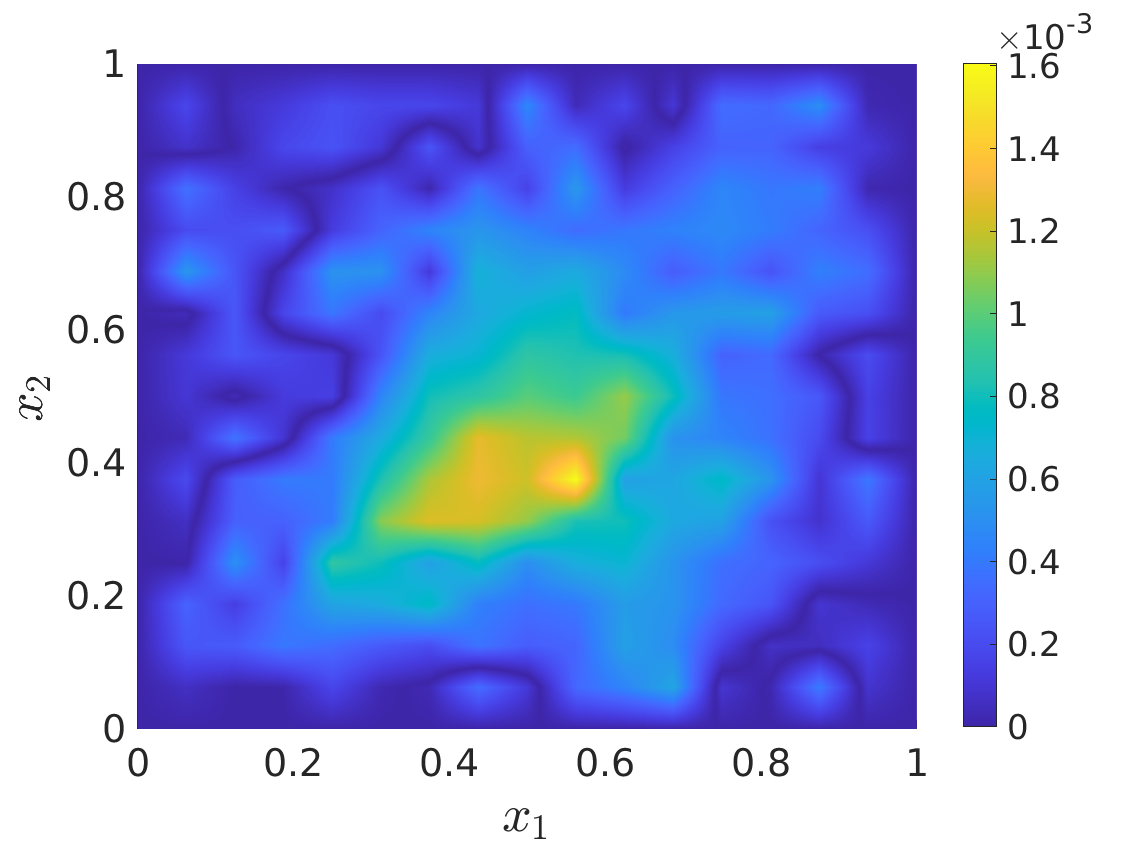}
     \end{subfigure} 
     \begin{subfigure}[a]{0.45\textwidth}
         \centering
         \includegraphics[width=\textwidth]{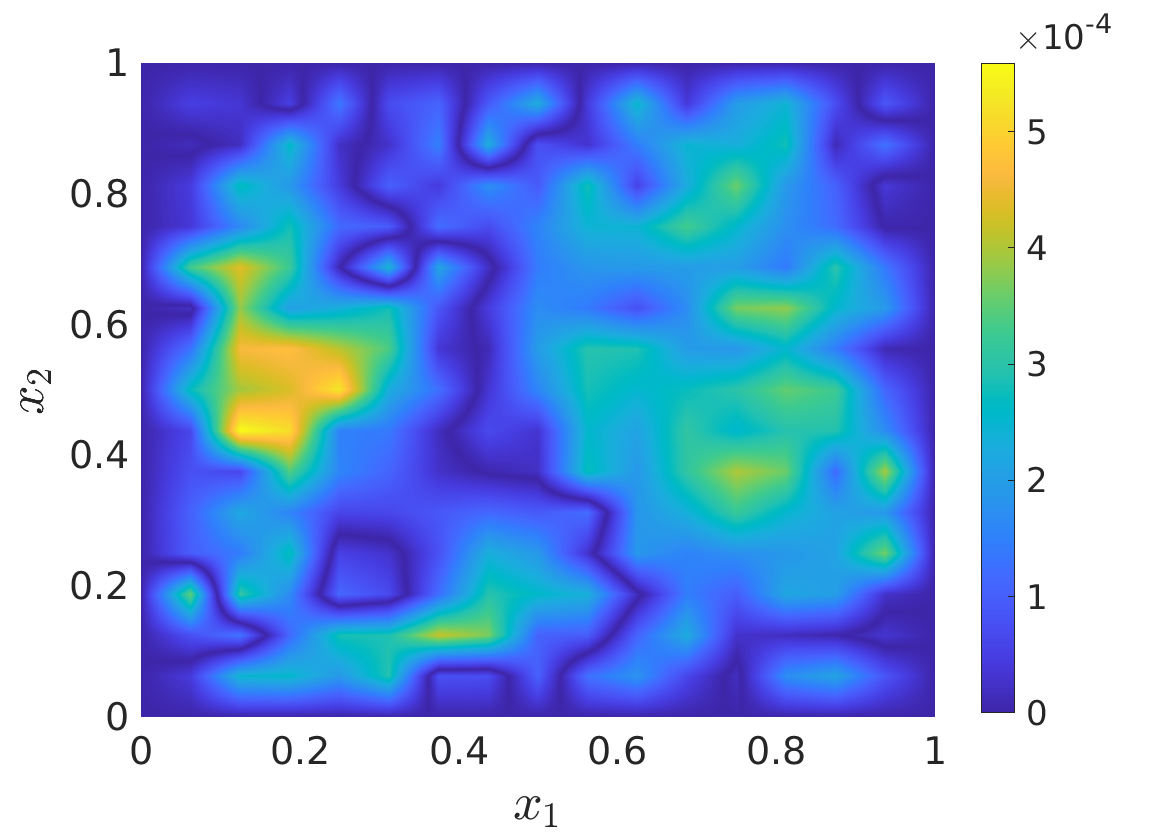}
     \end{subfigure} 
     \begin{subfigure}[a]{0.45\textwidth}
         \centering
         \includegraphics[width=\textwidth]{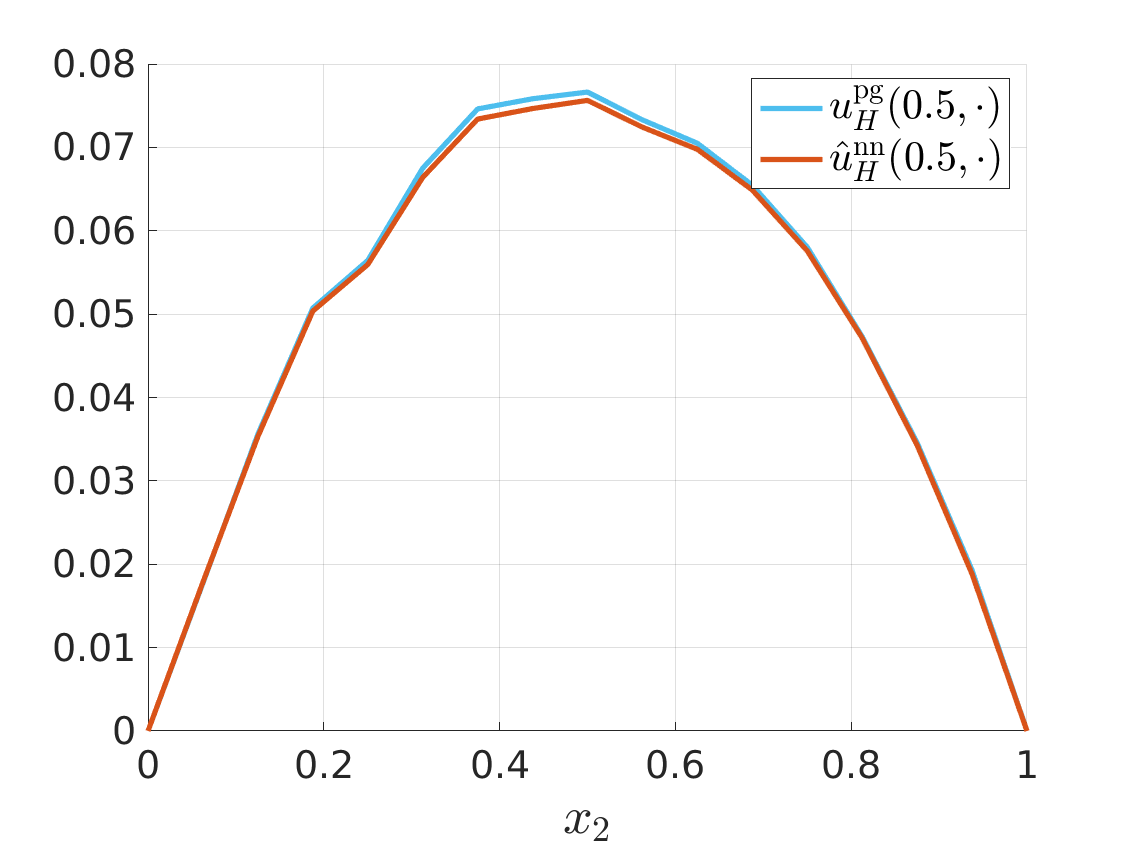}
     \end{subfigure}
     \begin{subfigure}[a]{0.45\textwidth}
         \centering
         \includegraphics[width=\textwidth]{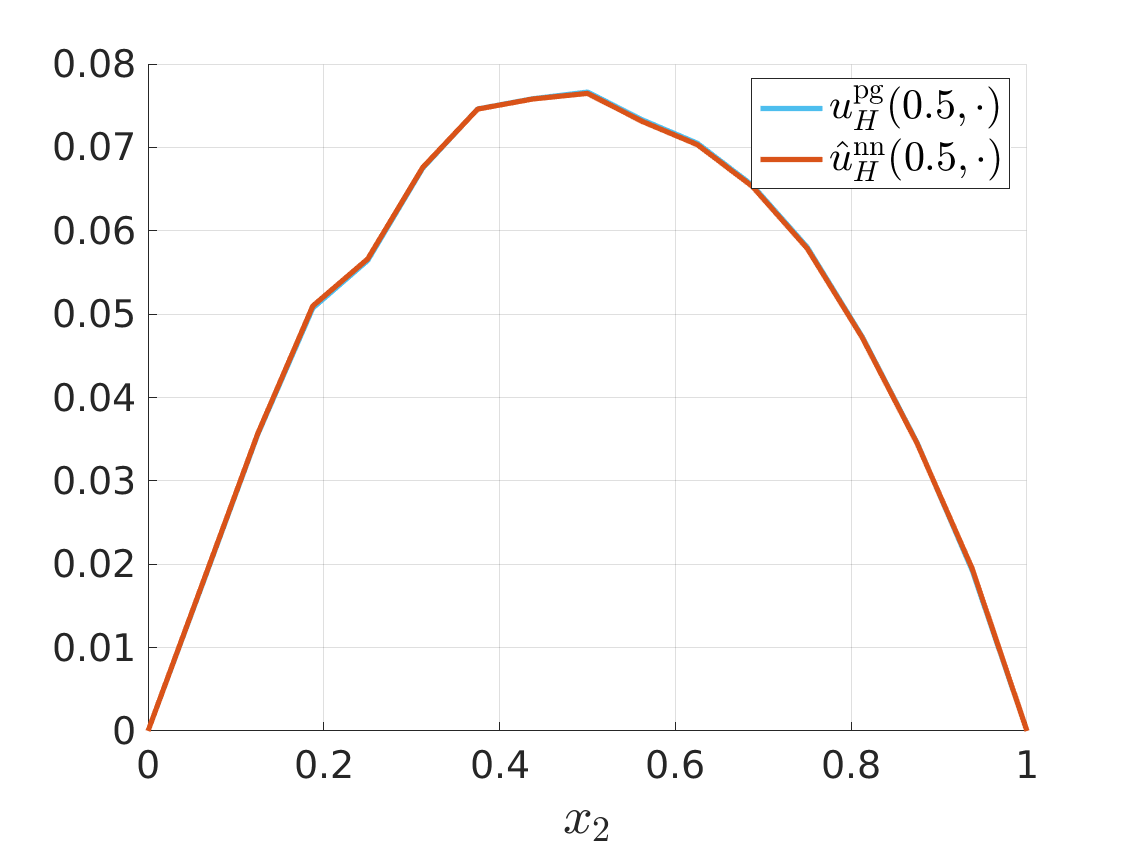}
     \end{subfigure}  
          \begin{subfigure}[a]{0.45\textwidth}
         \centering
         \includegraphics[width=\textwidth]{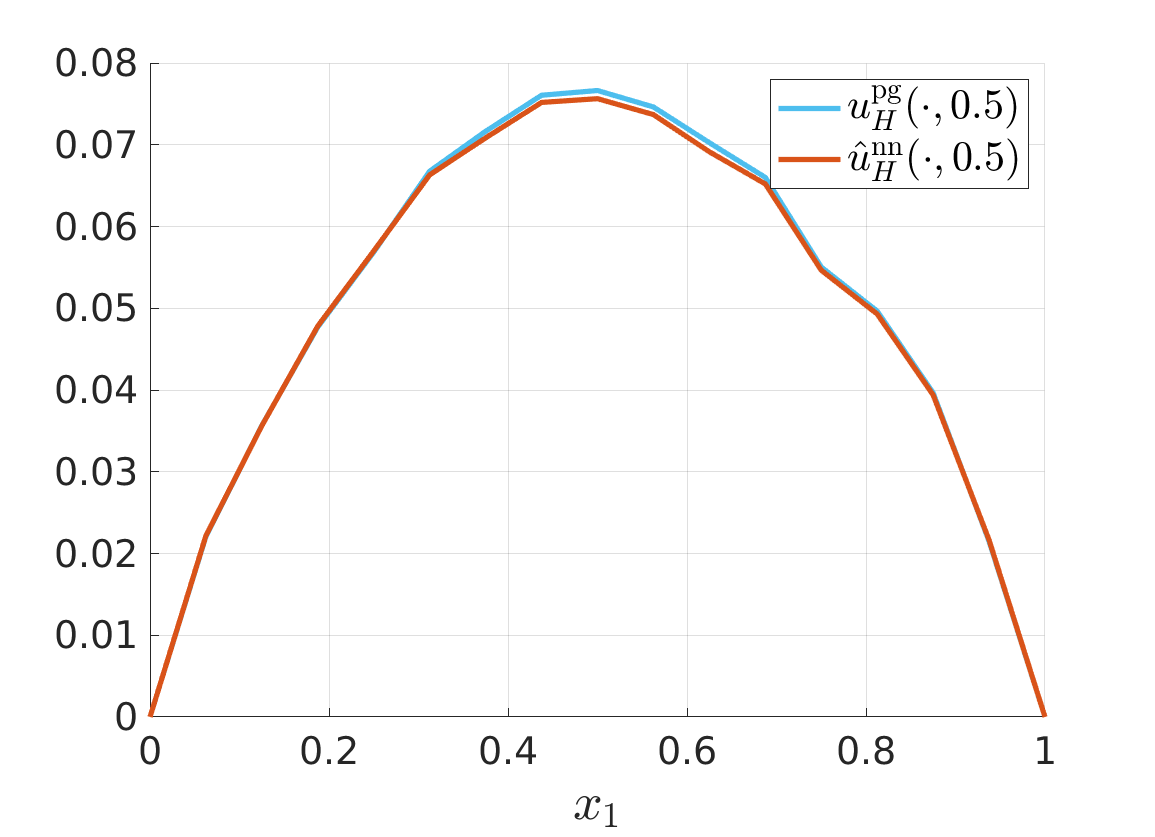}
     \end{subfigure}
     \begin{subfigure}[a]{0.45\textwidth}
         \centering
         \includegraphics[width=\textwidth]{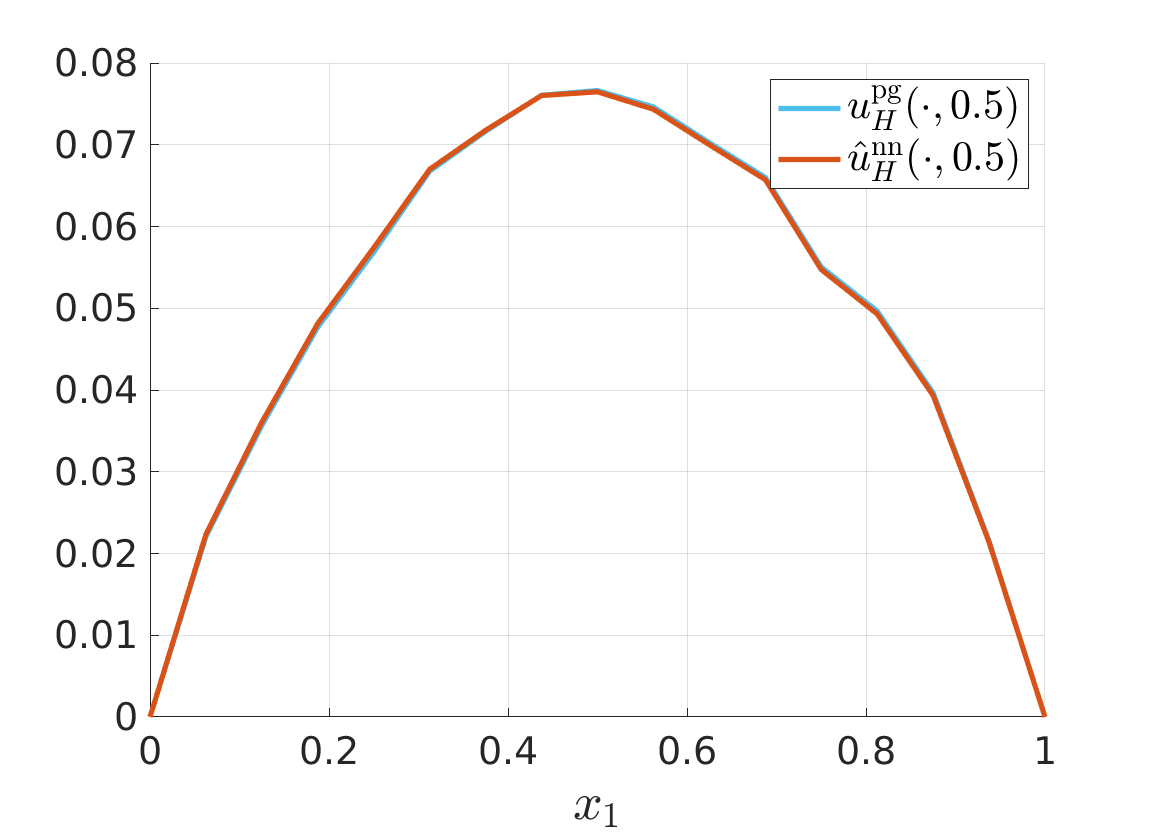}
     \end{subfigure}  
\caption{Experiment 1: The left column shows the results for the small dataset $\mathfrak{d}_{0.5}$: We see $\vert u_H^\mathrm{pg}(x) - \hat{u}_H^\mathrm{nn}(x) \vert$ (top row), the comparison of $u_H^\mathrm{pg}$ vs.~$\hat{u}_H^\mathrm{nn}$ along the cross-sections $x_1 = 0.5$ (middle row) and $x_2 = 0.5$ (bottom row). The right column shows the results for large dataset $\mathfrak{D}_{0.5}$.}
\label{fig:exploco}
\end{figure}

In this experiment, we retrain $\widehat{\Psi}^\mathrm{pg}$ on the low-contrast datasets $\mathfrak{d}_{0.5}$ and $\mathfrak{D}_{0.5}$.
In Figure~\ref{fig:lossloco}, we plot the development of the loss function $\mathcal{J}$ in \eqref{eq:loss} on the training and validation sets during the $60$ respective $40$ epochs of training.
We subsequently present the outcomes for both the small and large datasets.\\
For the \textbf{small dataset} $\mathfrak{d}_{0.5}$, we achieve a training loss of approximately $1.47 \cdot 10^{-4}$ and a validation loss of roughly $6.77\cdot 10^{-4}$ after $60$ epochs. 
How do these results compare to the uniformly elliptic case studied in~\cite{KroMP21}?
On the training set, we obtain an average relative error that is about half an order of magnitude larger; on the validation set, the average relative error is one order of magnitude larger compared to the result in~\cite{KroMP21}. This is also confirmed for the loss on the test set, which is approximately $6.88 \cdot 10^{-4}$. 
To test the network's applicability for our model problem~\eqref{eq:stochweak}, we sample a new low-contrast realization $\mathcal{A}_{0.5}^{(t)}$ with a contrast of approximately $3.32 \cdot 10^3$. This is a representative example to test the network $\widehat{\Psi}_{0.5}^{\mathrm{pg}}$ since the contrast of most realizations of $\mathcal{A}_{0.5}$ is on the order of $\mathcal{O}(10^3)$. After computing the two surrogates based on the PG-LOD and the network $\widehat{\Psi}^{\mathrm{pg}}_{0.5}$, we compare the corresponding approximate solutions $u_H^{\mathrm{pg}}$ and $\hat{u}_H^\mathrm{nn}$ in Figure~\ref{fig:exploco}. We obtain an $L^2$ error $\|u_H^\mathrm{pg} - \hat{u}_H^\mathrm{nn} \|_{\LL} \approx 4.49 \cdot 10^{-4}$ and a spectral norm difference $\|(\mathbf{S}_{\mathcal{A}_{0.5}}^{\mathrm{pg}})^{(t)} - (\widehat{\mathbf{S}}_{\mathcal{A}_{0.5}}^{\mathrm{nn}})^{(t)}\|_2 \approx 2.45 \cdot 10^{-1}$. 

Hence, we expect that the neural-network-based approximate solution is close to the PG-LOD approximation, our ground truth.
The results are shown in the left column of Figure~\ref{fig:exploco}.
A visual inspection of the two approximate solutions along the cross-sections $x_1 = 0.5$ and $x_2 = 0.5$ shows that, as expected, the neural network-based approximation is close to the ground truth. However, we also observe a systematic error in the neural network approximation, which undershoots the reference PG-LOD approximation.

For the \textbf{large dataset} $\mathfrak{D}_{0.5}$, we obtain a training loss of $1.84 \cdot 10^{-4}$ and a validation loss of $4.29 \cdot 10^{-4}$ after $40$ epochs. 
The difference between training and validation loss is smaller compared to the dataset $\mathfrak{d}_{0.5}$, which indicates less overfitting. In addition, we obtain a smaller test loss of $4.34 \cdot 10^{-4}$.  
This time, we obtain an $L^2$ error $\|u_H^\mathrm{pg} - \hat{u}_H^\mathrm{nn} \|_{\LL} \approx 1.18 \cdot 10^{-4}$ and a spectral norm difference $\|(\mathbf{S}_{\mathcal{A}_{0.5}}^{\mathrm{pg}})^{(t)} - (\widehat{\mathbf{S}}_{\mathcal{A}_{0.5}}^{\mathrm{nn}})^{(t)}\|_2 \approx 1.50 \cdot 10^{-1}$. 

Thus, we expect better results for new realizations. 
The results are shown in the right column of Figure~\ref{fig:exploco} and are consistent with our expectations: we observe a visible improvement compared to the small dataset. In particular, the systematic error is almost completely eliminated.

We conclude that the transfer learning approach yields good results for lognormal coefficients with a low contrast in the order of $\mathcal{O}(10^3)$, provided that the dataset is large enough.
Compared to the uniformly elliptic setting in \cite{KroMP21}, the results are slightly worse; however, this is expected, as the lognormal diffusion is much more challenging for numerical solvers. 

In the next experiment, we increase the contrast of the underlying coefficient class and investigate the impact on the neural network surrogate.

\subsection{Experiment 2: Moderate-contrast realizations}\label{subsec:mico}
In this experiment, we retrain the network $\widehat{\Psi}^{\mathrm{pg}}$ on the moderate-contrast datasets $\mathfrak{d}_1$ and $\mathfrak{D}_1$. 

\begin{figure}
     \centering
     \begin{subfigure}[a]{0.49\textwidth}
         \centering
         \includegraphics[width=\textwidth,trim= 50 0 100 0,clip]{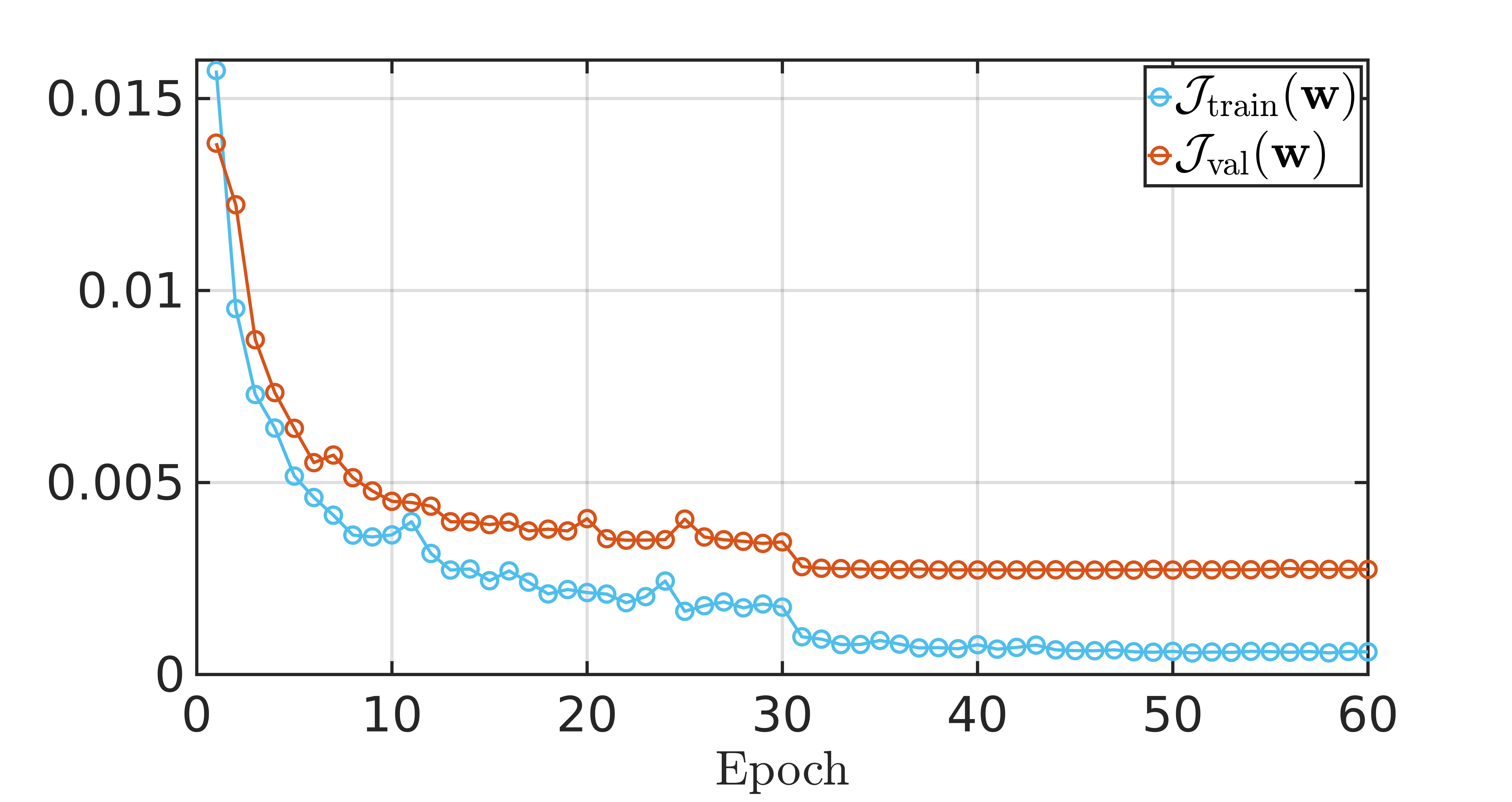}
     \end{subfigure} 
          \begin{subfigure}[a]{0.49\textwidth}
         \centering
         \includegraphics[width=\textwidth,trim= 50 0 100 0,clip]{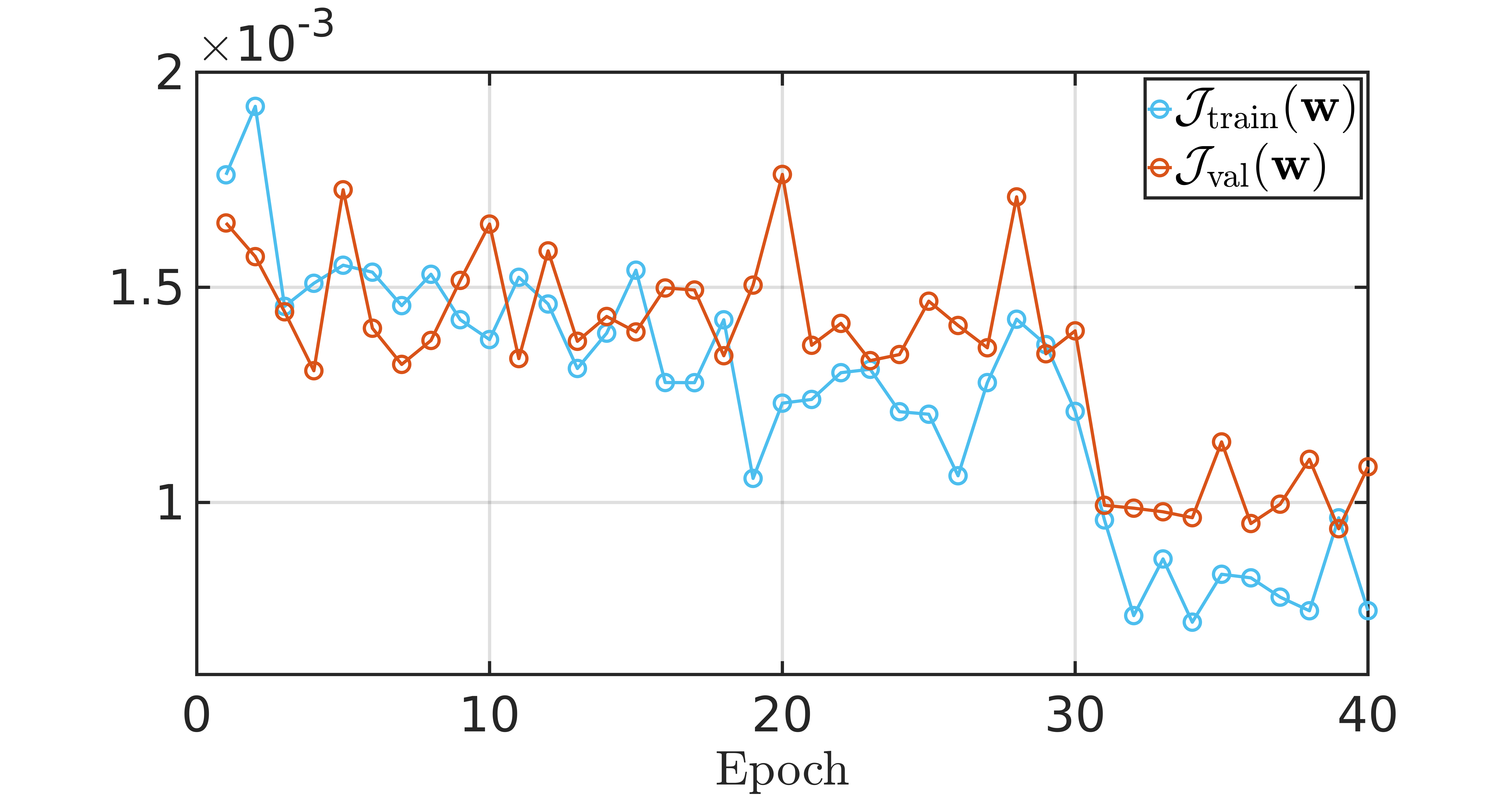}
     \end{subfigure} 
\caption[Development of loss for moderate-contrast dataset.]{Experiment 2: Training and validation loss for the moderate-contrast datasets~$\mathfrak{d}_{1}$ (left) and~$\mathfrak{D}_{1}$ (right) .}
\label{fig:lossmico}
\end{figure}
\begin{figure}
     \centering
     \begin{subfigure}[a]{0.49\textwidth}
         \centering
         \includegraphics[width=\textwidth]{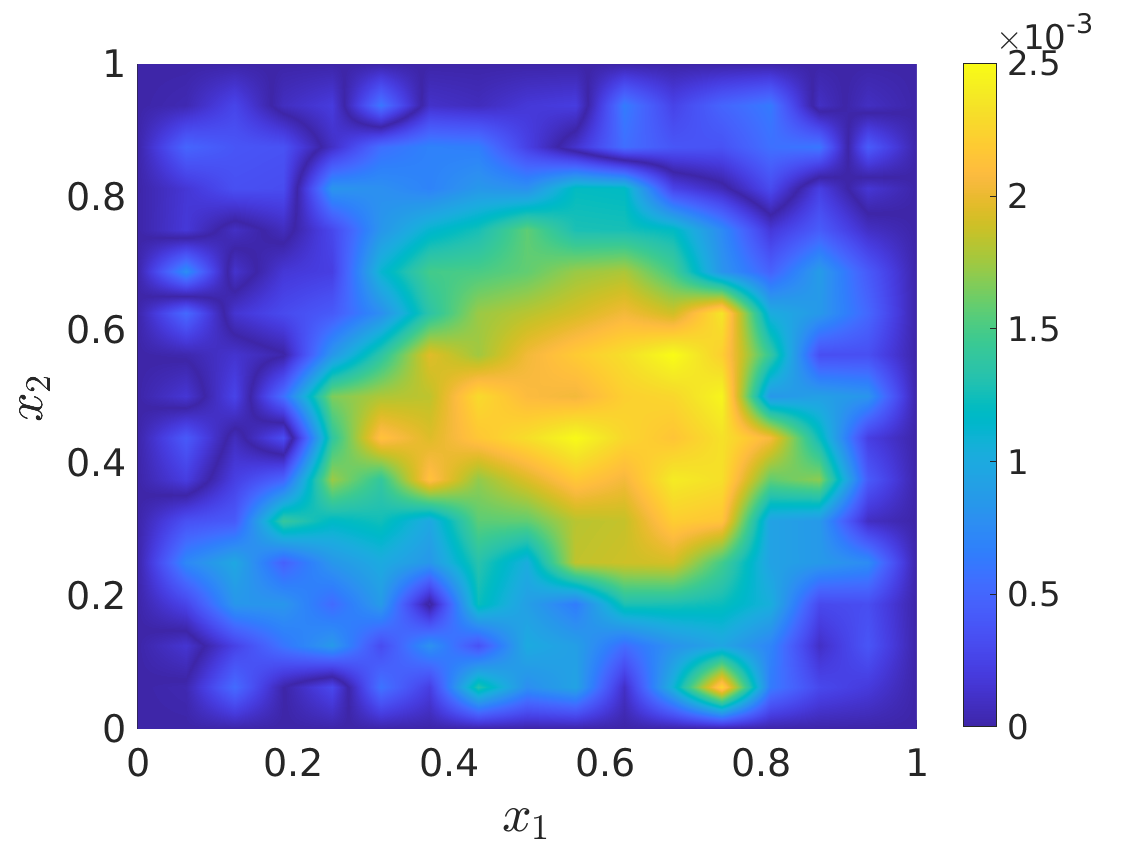}
     \end{subfigure} 
     \begin{subfigure}[a]{0.49\textwidth}
         \centering
         \includegraphics[width=\textwidth]{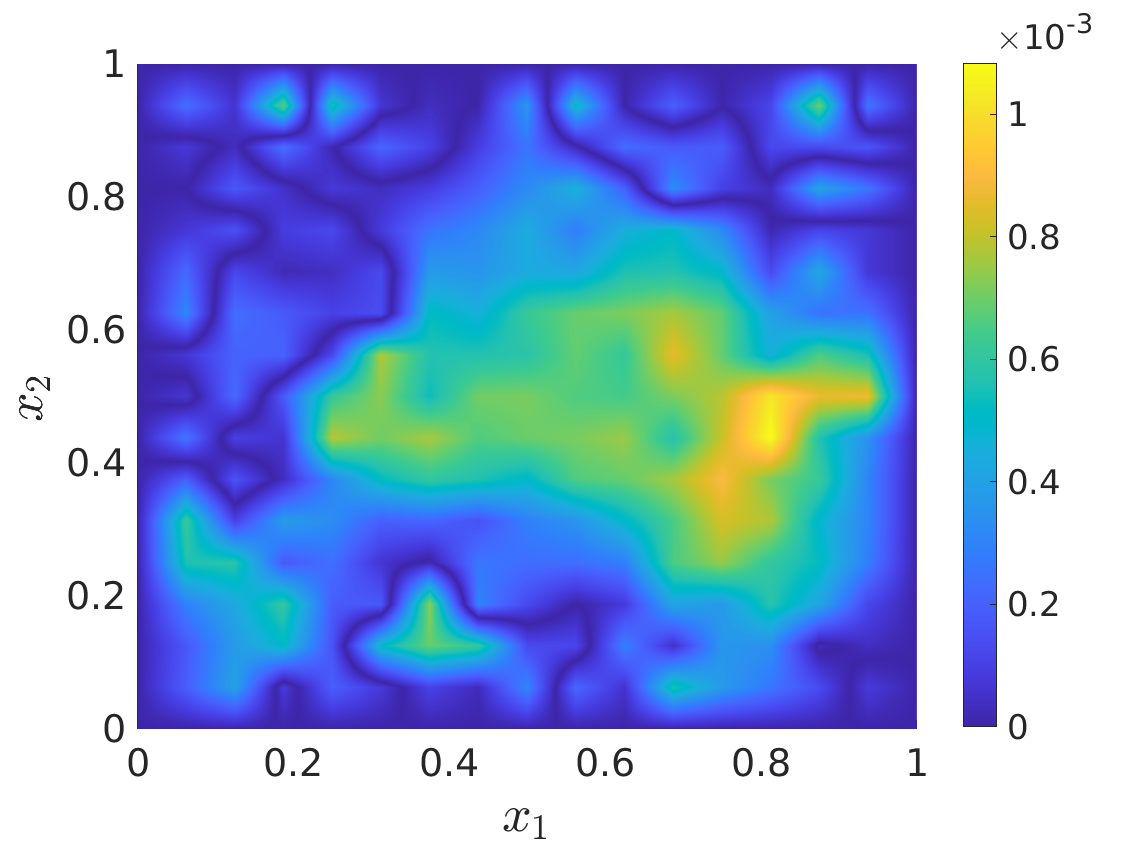}
     \end{subfigure} 
     \begin{subfigure}[a]{0.49\textwidth}
         \centering
         \includegraphics[width=\textwidth]{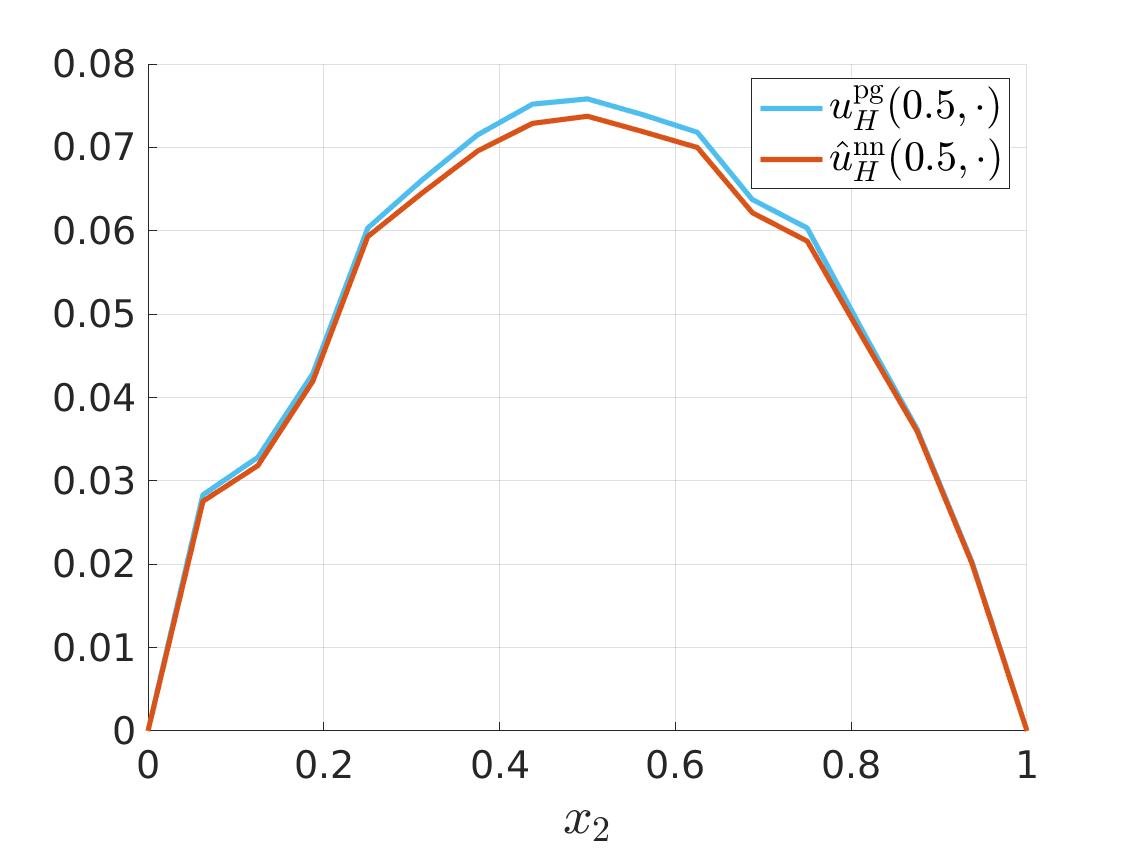}
     \end{subfigure}
     \begin{subfigure}[a]{0.49\textwidth}
         \centering
         \includegraphics[width=\textwidth]{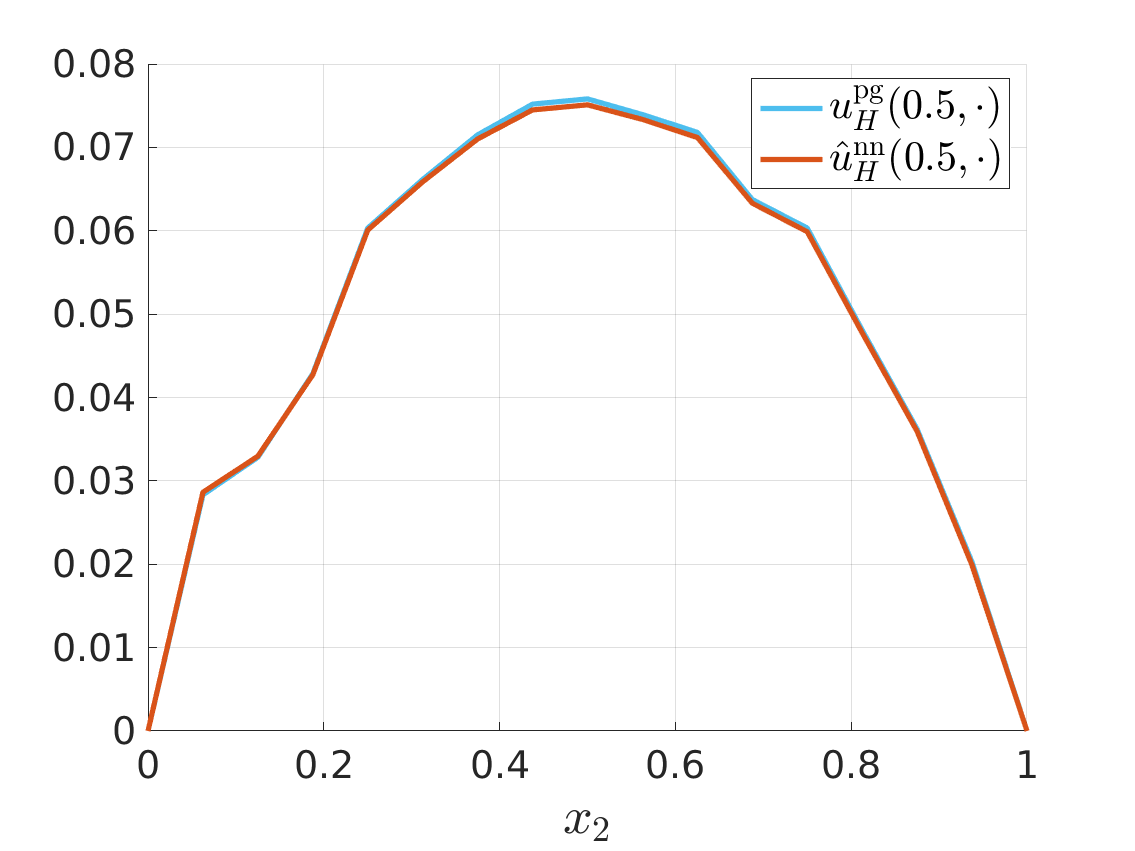}
     \end{subfigure}  
          \begin{subfigure}[a]{0.49\textwidth}
         \centering
         \includegraphics[width=\textwidth]{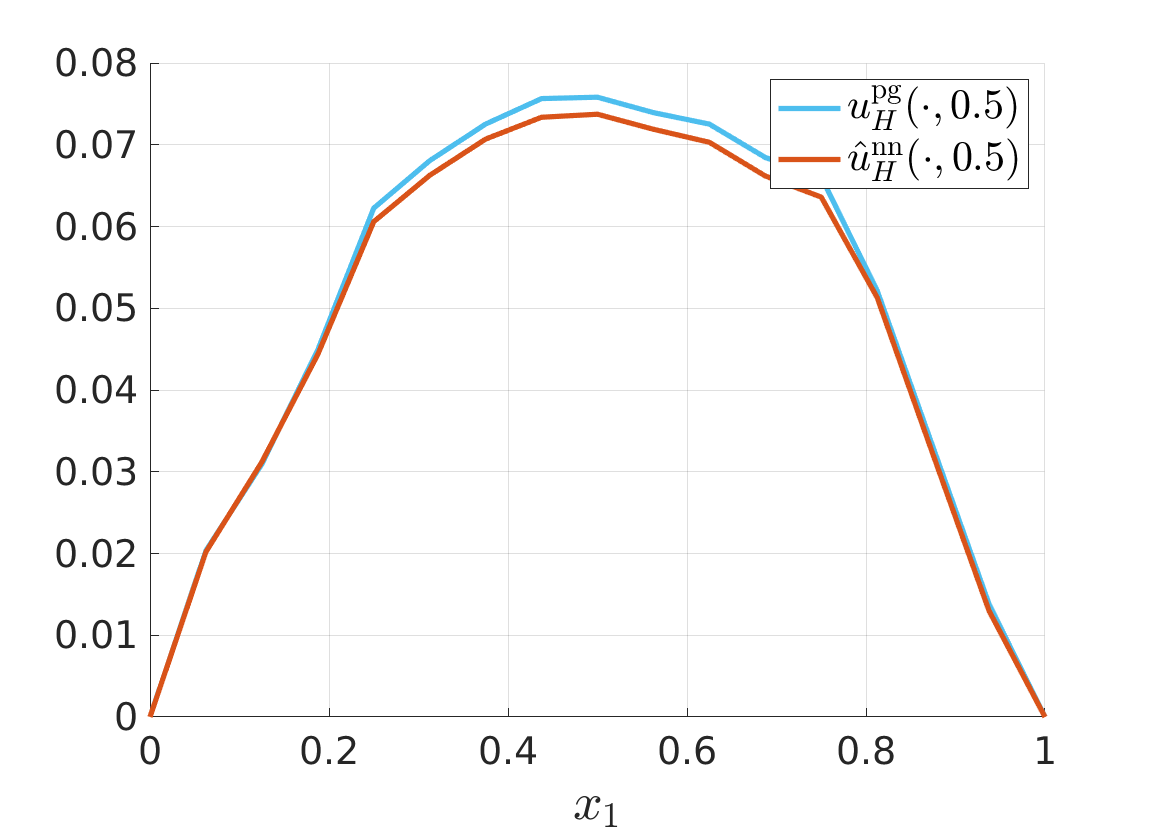}
     \end{subfigure}
     \begin{subfigure}[a]{0.49\textwidth}
         \centering
         \includegraphics[width=\textwidth]{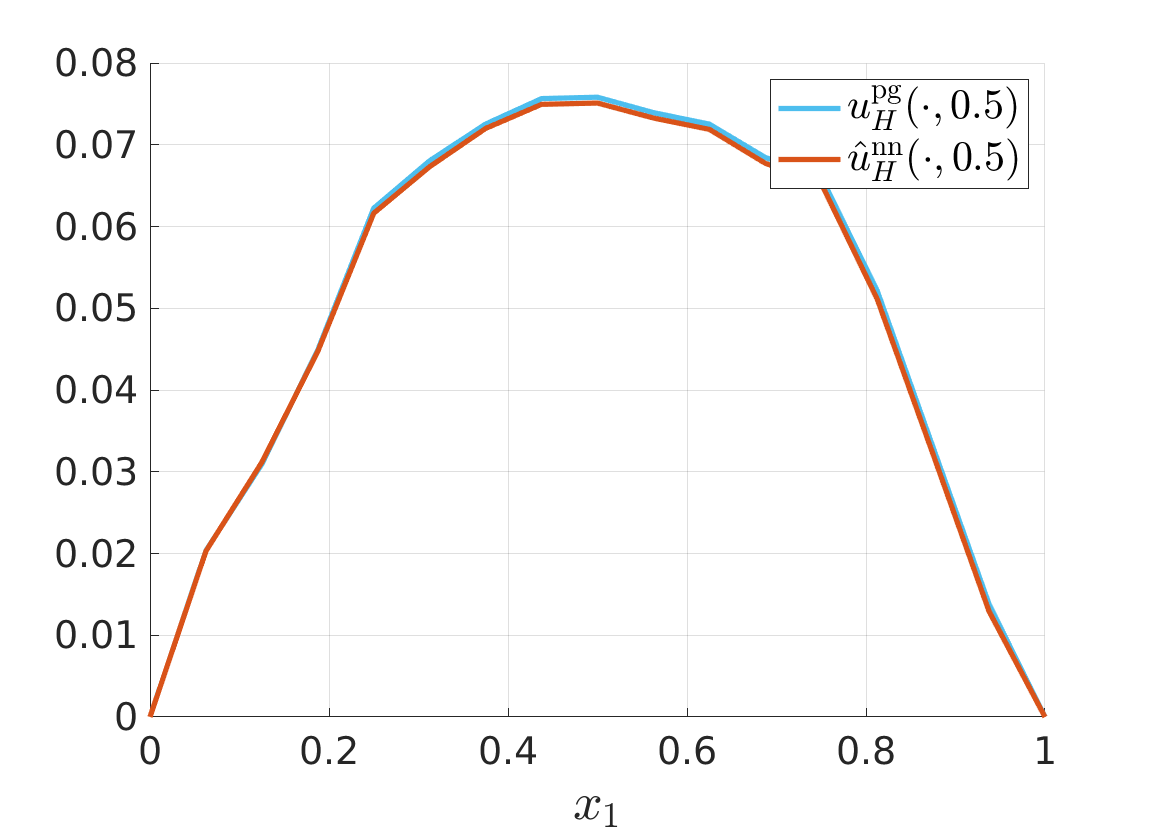}
     \end{subfigure}  
\caption{Experiment 2: The left column shows the results for the small dataset $\mathfrak{d}_1$, namely the error $\vert u_H^\mathrm{pg}(x) - \hat{u}_H^\mathrm{nn}(x) \vert$ (top row), the comparison of $u_H^\mathrm{pg}$ vs.~$\hat{u}_H^\mathrm{nn}$ along the cross-sections $x_1 = 0.5$ (middle row) and $x_2 = 0.5$ (bottom row). The right column shows the results for large dataset $\mathfrak{D}_1$.}
\label{fig:expmico}
\end{figure}

We begin by discussing the results for the \textbf{small dataset}  $\mathfrak{d}_1$.
We first examine the development of the loss function $\mathcal{J}$ in \eqref{eq:loss}, which is shown in Figure~\ref{fig:lossmico} on the left-hand side. After $60$ epochs, we arrive at a training loss of approximately $5.88 \cdot 10^{-4}$ and a validation loss of about $2.74 \cdot 10^{-3}$.
In both cases, this is about half an order of magnitude larger than in the low-contrast setting of Experiment 1. Moreover, the difference between the final training and validation losses is larger than in the previous setting. This suggests that, with increasing contrast, the network's ability to fit a given training dataset decreases.
In addition, the capability of the network to generalize beyond the training data decreases. 
The loss on the test set $\mathfrak{d}_1^{\mathrm{test}}$ is again close to the validation loss, with an approximate value of $2.81\cdot 10^{-3}$. 
To test the practical performance of the retrained network $\widehat{\Psi}^{\mathrm{pg}}_1$, we sample a new realization $\mathcal{A}_{1}^{(t)}$ of the random field $\mathcal{A}_{1}$, which has a contrast of approximately $2.25 \cdot 10^4$.
This is one order of magnitude larger than the test realization in Experiment 1. Based on this realization, we obtain an $L^2$ error $\|u_H^\mathrm{pg} - \hat{u}_H^\mathrm{nn} \|_{\LL} \approx 1.08 \cdot 10^{-3}$ and a spectral norm difference $\|(\mathbf{S}_{\mathcal{A}_{1}}^{\mathrm{pg}})^{(t)} - (\widehat{\mathbf{S}}_{\mathcal{A}_{1}}^{\mathrm{nn}})^{(t)}\|_2 \approx 4.75 \cdot 10^{-1}$ between the surrogates corresponding to the PG-LOD and $\widehat{\Psi}^{\mathrm{pg}}_1$, respectively. We visualize the results of this experiment in the left column of Figure~\ref{fig:expmico}. We clearly see that the difference between the neural-network-based approximation and the PG-LOD-based  ground truth is more pronounced. The qualitative behavior of the approximation is correct; however, the neural-network-based approximation undershoots the reference solution in terms of the amplitude along both cross-sections. We also performed test runs with many other realizations of $\mathcal{A}_1$ and observed the same qualitative behavior: the amplitude of the neural-network-based approximate solution was consistently either too large or too small. 

We see in the right column of Figure~\ref{fig:expmico} that the mismatch between the neural-network-based solution and the PG-LOD-based ground truth can be decreased by using a significantly \textbf{larger dataset}. 
Note that the dataset $\mathfrak{D}_1$ contains more than 10 times more realizations compared to the dataset $\mathfrak{d}_1$. Again, the training and validation loss are much closer during training, indicating less overfitting. The final training loss is $7.48\cdot 10^{-4}$; the final validation loss is $1.1 \cdot 10^{-3}$. The error on the test set after $40$ epochs amounts to $7.27 \cdot 10^{-4}$. For the test coefficient $\mathcal{A}_1^{(t)}$, the spectral error between the global PG-LOD and the NN-based stiffness matrix is $0.56$, and the L2 error between the corresponding PDE solutions is $3.81\cdot 10^{-4}$.

We conclude that for medium-contrast lognormal diffusion coefficients, we require a significantly larger dataset for the training of the neural network surrogate to obtain results with a similar accuracy as for low-contrast lognormal coefficients.

\subsection{Experiment 3: High-contrast realizations}\label{subsec:hico}
\begin{figure}
     \centering
     \begin{subfigure}[a]{0.495\textwidth}
         \centering
         \includegraphics[width=\textwidth,trim= 80 0 100 0,clip]{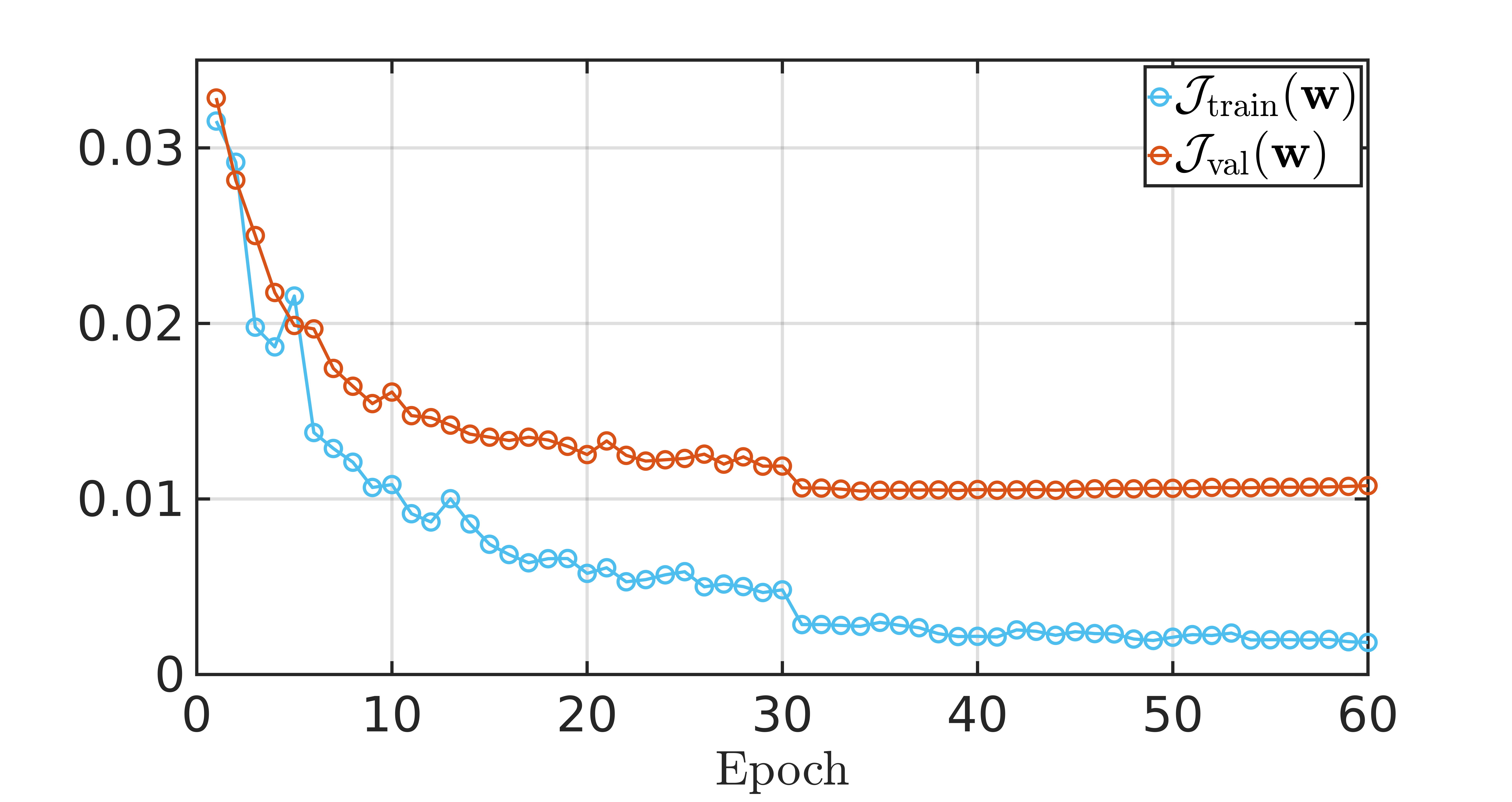}
     \end{subfigure} 
    \begin{subfigure}[a]{0.495\textwidth}
         \centering
         \includegraphics[width=\textwidth,trim= 80 0 100 0,clip]{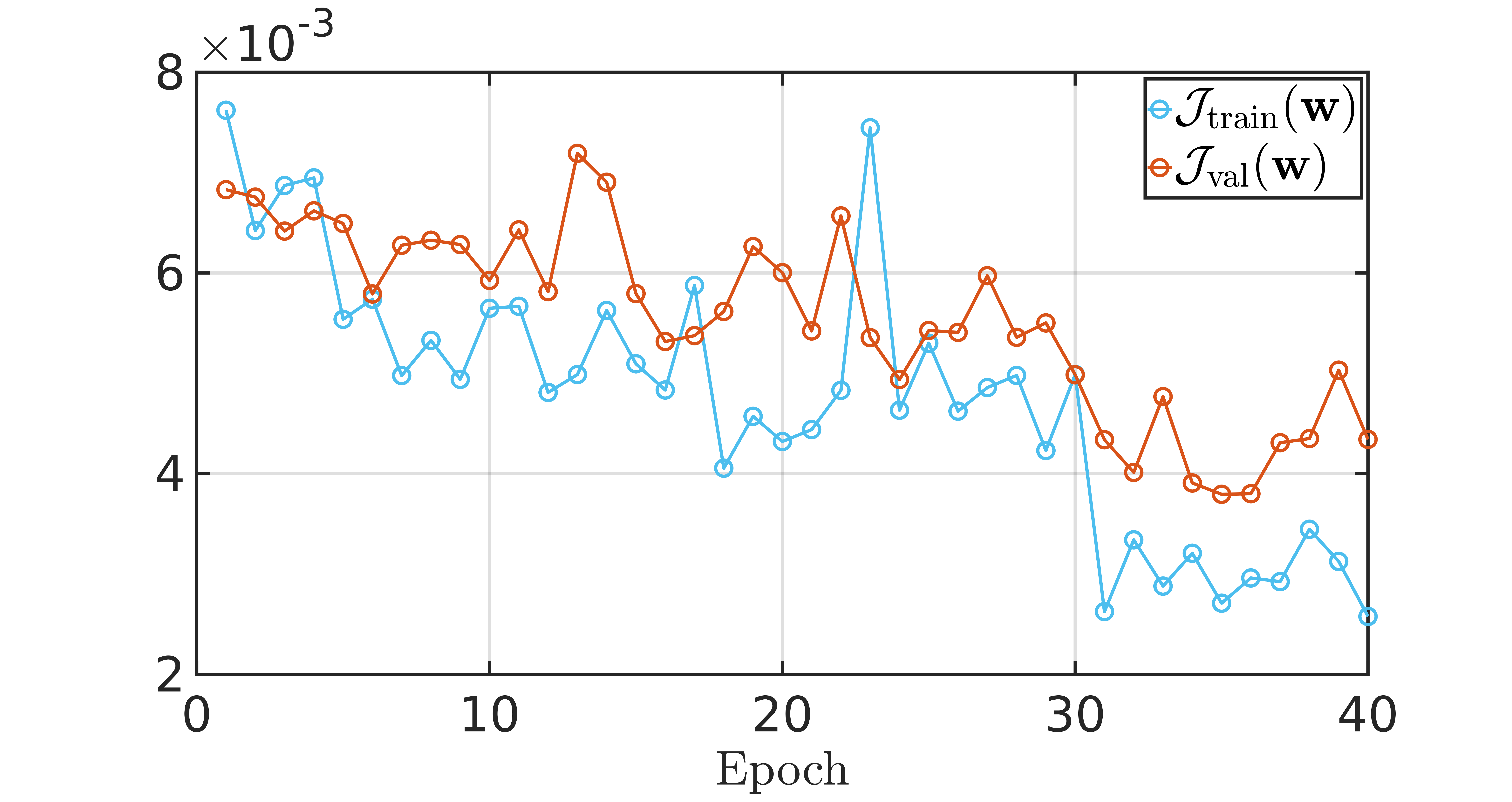}
     \end{subfigure} 
\caption[Development of loss for high-contrast dataset.]{Experiment 3: Training and validation loss for the high-contrast datasets~$\mathfrak{d}_{2}$ (left) and~$\mathfrak{D}_{2}$ (right) .}
\label{fig:losshico}
\end{figure}

\begin{figure}
     \centering
     \begin{subfigure}[a]{0.45\textwidth}
         \centering
         \includegraphics[width=\textwidth]{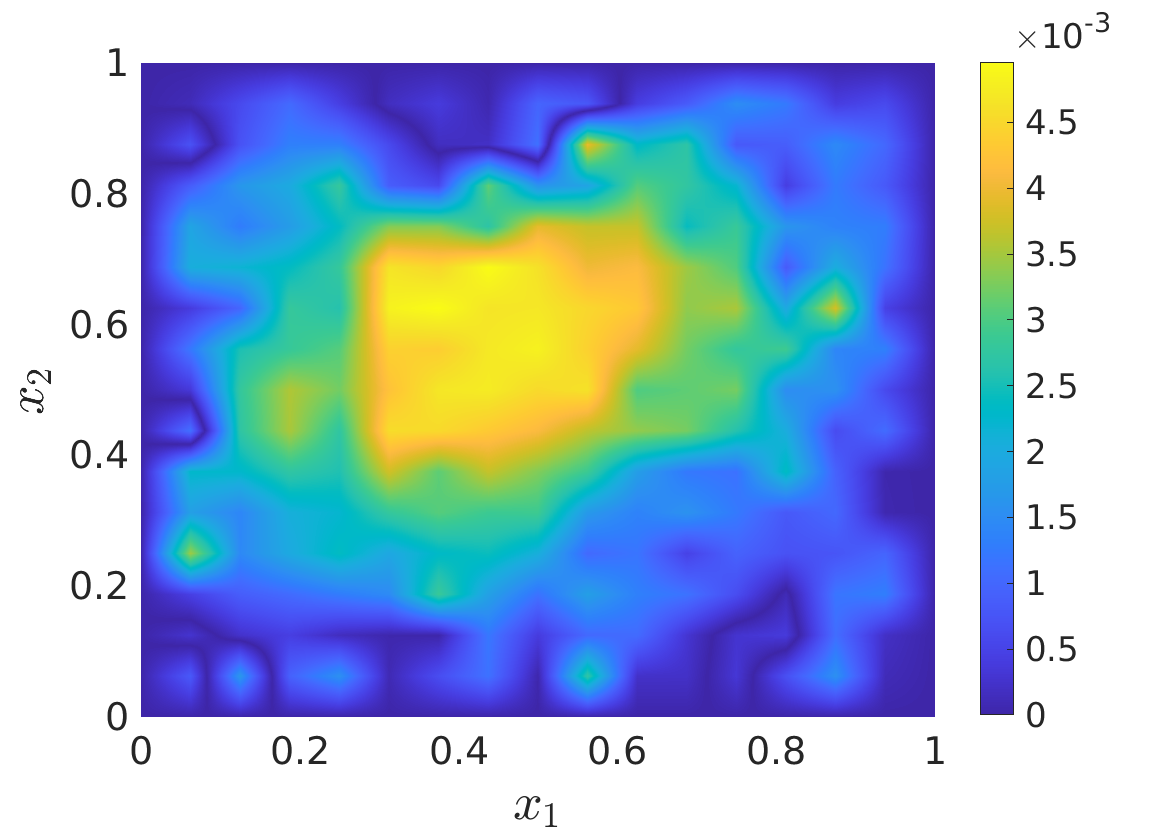}
     \end{subfigure} 
     \begin{subfigure}[a]{0.45\textwidth}
         \centering
         \includegraphics[width=\textwidth]{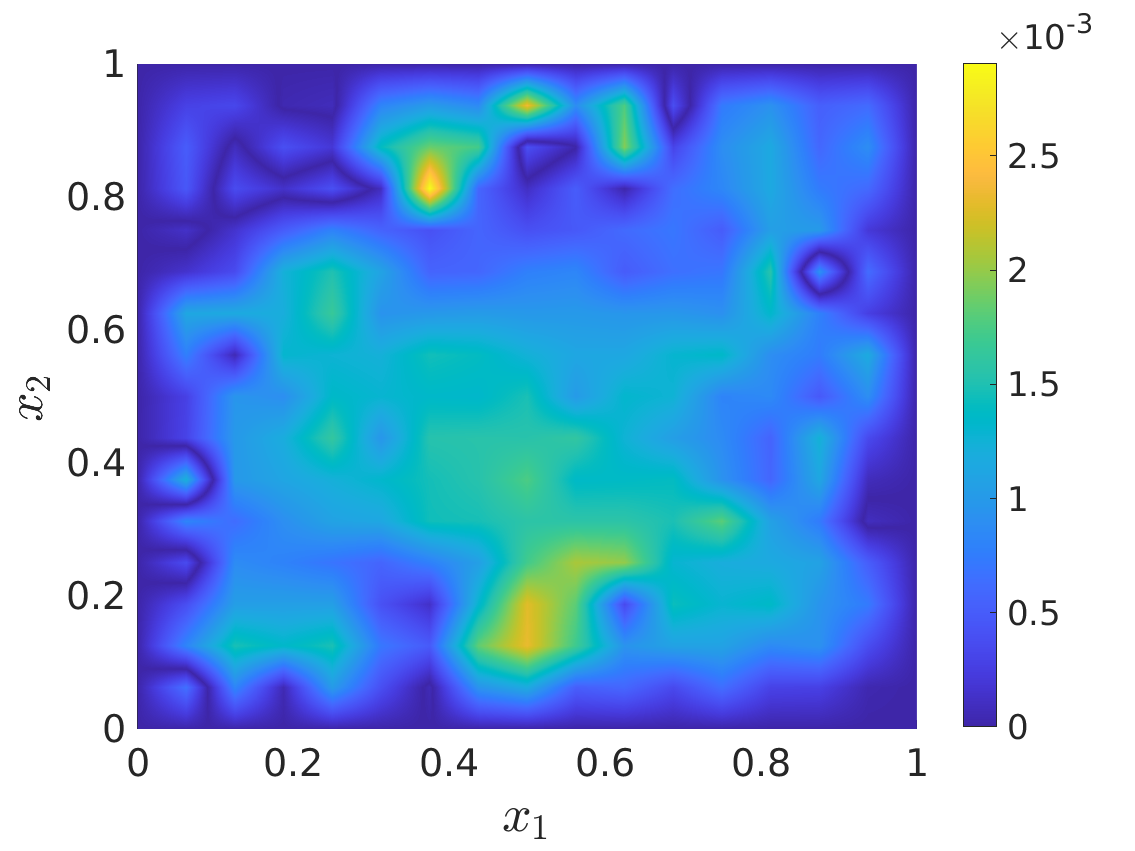}
     \end{subfigure} 
     \begin{subfigure}[a]{0.45\textwidth}
         \centering
         \includegraphics[width=\textwidth]{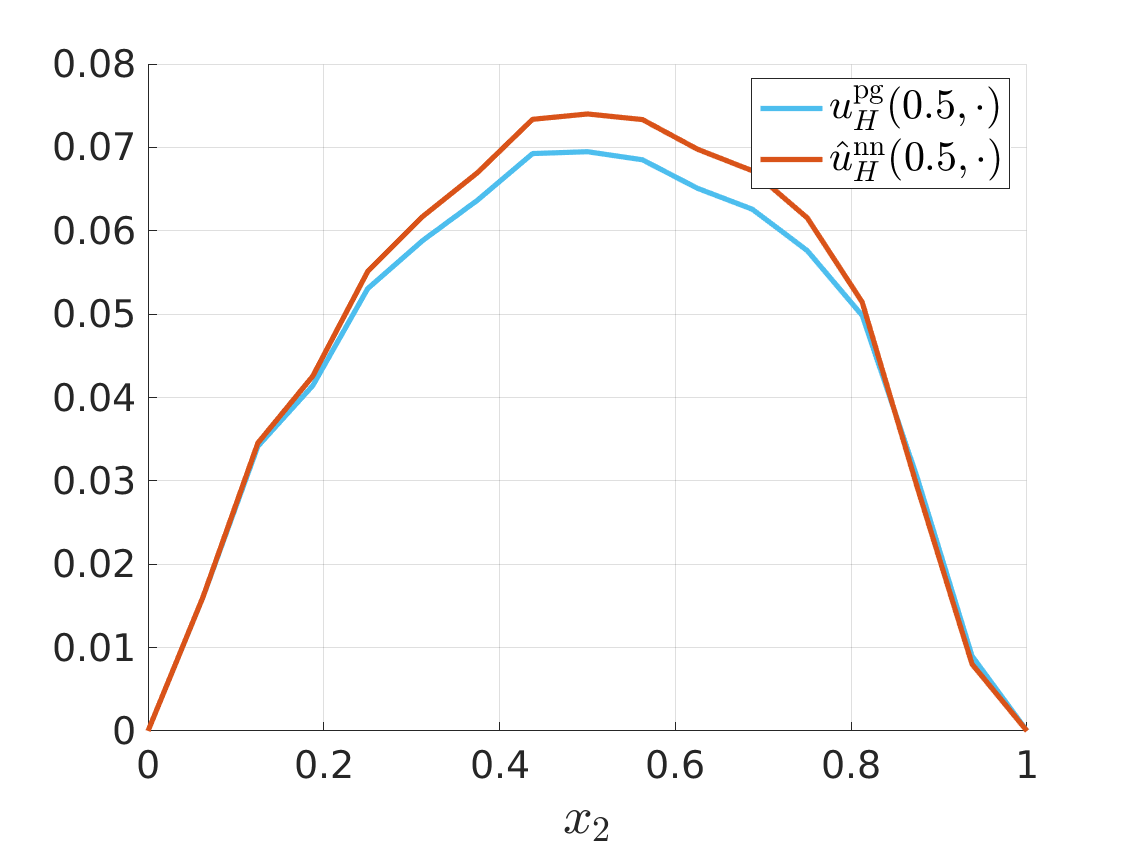}
     \end{subfigure}
     \begin{subfigure}[a]{0.45\textwidth}
         \centering
         \includegraphics[width=\textwidth]{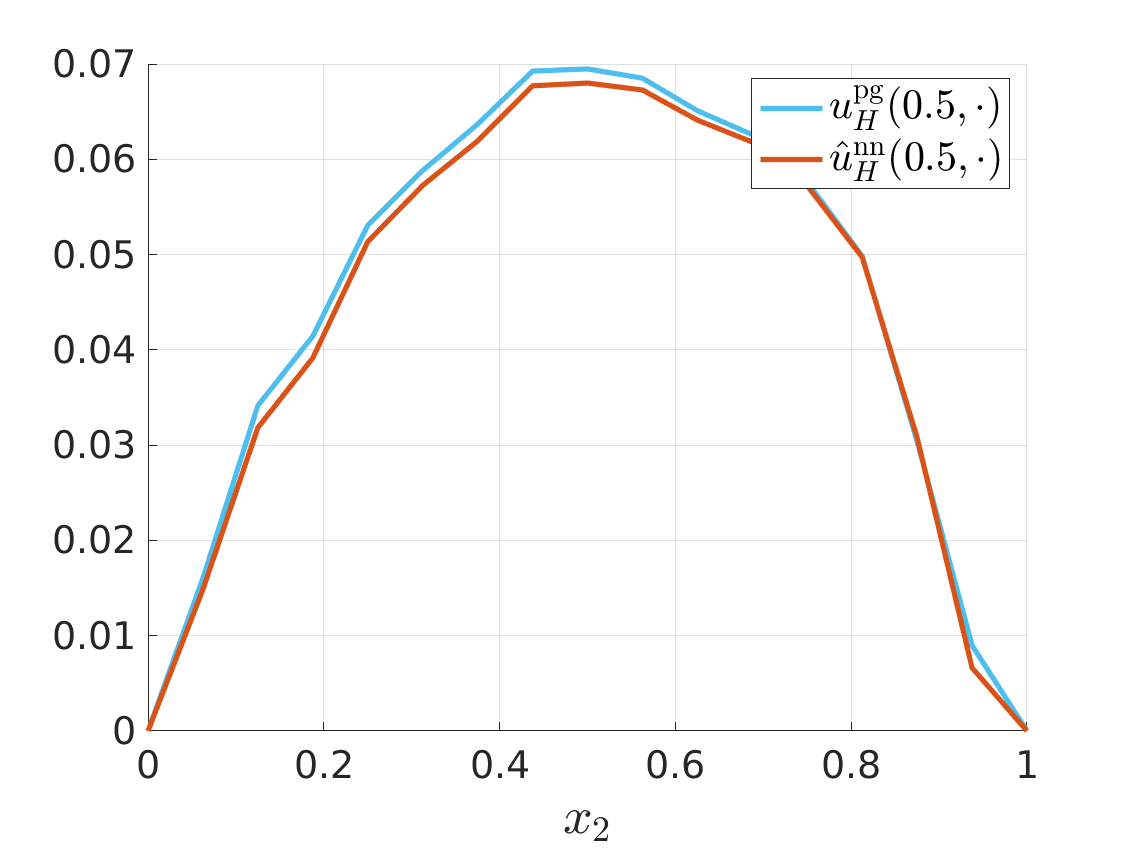}
     \end{subfigure}  
          \begin{subfigure}[a]{0.45\textwidth}
         \centering
         \includegraphics[width=\textwidth]{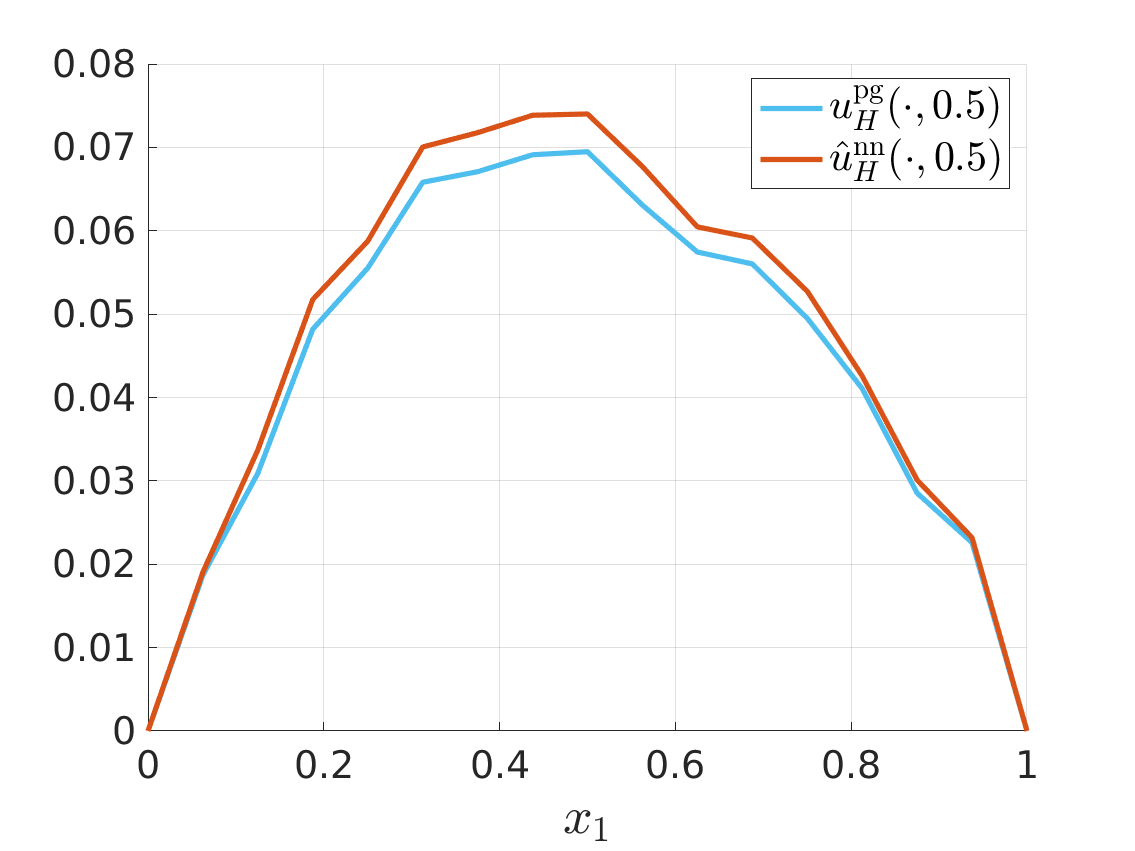}
     \end{subfigure}
     \begin{subfigure}[a]{0.45\textwidth}
         \centering
         \includegraphics[width=\textwidth]{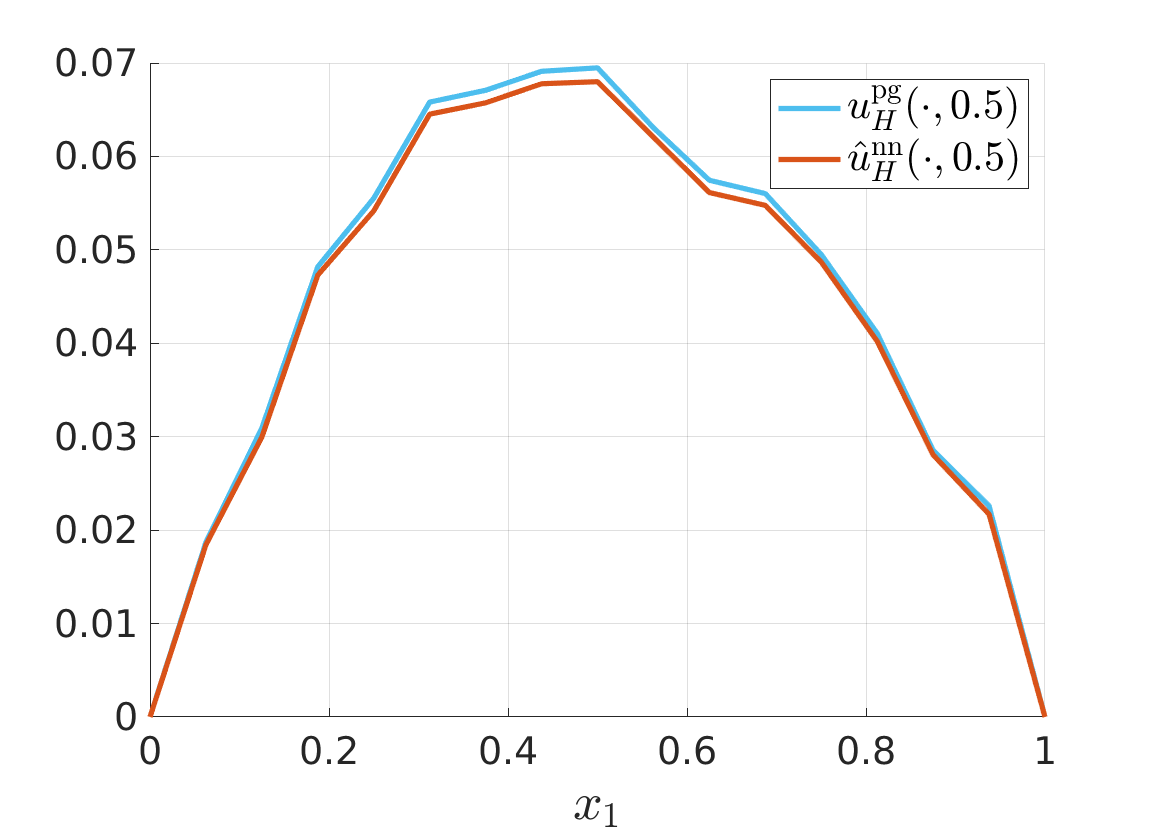}
     \end{subfigure}  
\caption{Experiment 3: The left column shows the results for the small dataset $\mathfrak{d}_2$, namely the error $\vert u_H^\mathrm{pg}(x) - \hat{u}_H^\mathrm{nn}(x) \vert$ (top row), the comparison of $u_H^\mathrm{pg}$ vs.~$\hat{u}_H^\mathrm{nn}$ along the cross-sections $x_1 = 0.5$ (middle row) and $x_2 = 0.5$ (bottom row). The right column shows the results for large dataset $\mathfrak{D}_2$.}
\label{fig:exphico}
\end{figure}

Finally, we retrain the network $\widehat{\Psi}^{\mathrm{pg}}$ on the challenging datasets $\mathfrak{d}_2$ and $\mathfrak{D}_2$ with high-contrast realizations.
Figure~\ref{fig:losshico} shows the loss function $\mathcal{J}$ in \eqref{eq:loss}. For the small dataset, we again see that the training and validation losses drift further apart with increased contrast. After epoch $30$, we also observe a tendency for overfitting: the validation loss remains constant while the training loss continues to decrease. 
We begin by discussing the results for the \textbf{small dataset} $\mathfrak{d}_2$.
After $60$ epochs, we arrive at a training loss of roughly $1.83 \cdot 10^{-3}$ and a validation loss of approximately $1.08 \cdot 10^{-2}$. The loss on the test set in this setting amounts to $1.1\cdot 10^{-2}$. This means that for all three subsets of the dataset $\mathfrak{d}_2$, the average relative error is again half an order to an order of magnitude larger than in the previous experiment. Again, we generate an additional realization $\mathcal{A}_{2}^{(t)}$ for testing purposes, which, in this case, has an approximate contrast of $5.63 \cdot 10^5$. The $L^2$ error is $\|u_H^\mathrm{pg} - \hat{u}_H^\mathrm{nn} \|_{\LL} \approx 2.2 \cdot 10^{-3}$ and the spectral norm difference is $\|(\mathbf{S}_{\mathcal{A}_{2}}^{\mathrm{pg}})^{(t)} - (\widehat{\mathbf{S}}_{\mathcal{A}_{2}}^{\mathrm{nn}})^{(t)}\|_2 \approx 1.07$. 
A graphical representation of this test case is shown in Figure~\ref{fig:exphico}.
Looking at the cross-sections of the approximate solutions, it is clear that the accuracy of the neural-network-based approximation has deteriorated significantly. We observe a significant overshoot compared to the reference approximation based on the PG-LOD. 

For the \textbf{large dataset} $\mathfrak{D}_2$, the final training loss is $2.6\cdot 10^{-3}$, and the final validation loss is $4.3\cdot 10^{-3}$. The loss on the test set after $40$ epochs is $2.81\cdot 10^{-3}$. For the fresh test coefficient $\mathcal{A}_2^{(t)}$, we obtain a spectral difference of $0.63$ and an $L_2$ error of $9.37\cdot 10^{-4}$ between the respective approximations.

To conclude, the numerical experiments in this section support the following hypothesis: if we increase the contrast of the lognormal diffusion coefficients, then the neural-network-based approximation is less accurate, provided that the architecture and the number of training samples are fixed. However, this effect can be mitigated by increasing the amount of training data.

\section{Numerical experiments with hierarchical diffusion coefficients}\label{sec:numexp:hg}

In this section, the diffusion coefficient $\mathcal{A}$ in \eqref{eq:modelstrong} is given by $\mathcal{A}=\exp(\widetilde Z)$, where $\widetilde Z$ is a centered, hierarchical Gaussian random field with a random correlation length $\kappa$.
This means that for each fixed realization of $\kappa$, the random field $\widetilde Z(\cdot \vert \kappa)$ is a Gaussian random field with mean zero and Whittle--Matérn covariance function in \eqref{c}.
Importantly, the random field $\widetilde Z$ is not necessarily Gaussian.
Note that in this setting the realizations of $\mathcal{A}$ may still have a relatively large contrast; this is the same numerical challenge as for the lognormal coefficients in Section~\ref{sec:numexp:log}. However, the correlation length of the random field now changes across the realizations.
This should pose a manageable challenge for the neural-network-based LOD compression, which can handle relatively small correlation lengths, as we have seen in the previous section.
As in Section~\ref{sec:numexp:log}, we study the effects of an increasing contrast in the realizations of $\mathcal{A}$.
The experimental setup and simulation procedures are as described in Section~\ref{sec:numexp:log} unless stated otherwise.

Motivated by the examples in~\cite[Section 5]{Kressner2020}, the Gaussian random field $\widetilde Z(\cdot \vert \kappa)$ has the Whittle--Matérn covariance function in \eqref{c} with the hyperparameters $\sigma^2 \in \{0.5,1,2\}$ and $\nu = 1$. 
The correlation length ${\kappa} \sim \text{Unif}[2^{-6},2^{-3}]$ is now a uniformly distributed random variable taking values in the interval $[2^{-6},2^{-3}]$.

Again, we do not train a neural network from scratch but employ transfer learning. However, instead of initializing our training algorithm with the weights and biases of the neural network $\widehat{\Psi}^{\mathrm{pg}}$ in Section~\ref{sec:setup}, we use the neural network models trained in Experiments 1--3 (Section~\ref{sec:numexp:log}) as our initial model.
This is a reasonable approach since the correlation length $\kappa=2^{-6}$ in those experiments is the lower bound for the correlation lengths of the random field $\widetilde Z(\cdot \vert \kappa)$.

To build our datasets, we generate global realizations of the corresponding hierarchical random field for representative correlation lengths $\kappa\in\{2^{-3},2^{-4},2^{-5},2^{-6}\}$.
We generate $200$ realizations each for the low-contrast case and $800$ realizations each for the medium and high-contrast cases. 
The splitting into training, validation, and test datasets is performed according to the percentages specified in Section~\ref{sec:numexp:log}. 
That is, for each $\kappa$, we assign the first $160$ (resp. $640$) samples to the training set, the next $20$ (resp. $80$) samples to the validation set, and the last $20$ (resp. $80$) samples to the test set.
We split all samples into local subsamples $Z^{(i)}_{\kappa,T}$, take the point-wise exponential, and compute the corresponding local PG-LOD system matrices $\mathbf{S}^{(i)}_{\kappa,T}$. 

To train the networks, we use a batch size of $100$; we train each model for $60$ epochs and use the Adam optimizer with a step size of $10^{-3}$ until epoch $40$, and $10^{-4}$ thereafter.

In the following numerical experiments, we compare approximate solutions of the problem \eqref{eq:modelstrong} generated by the finite element method (FEM), PG-LOD, and our neural-network-based surrogate (NN-LOD) for one fresh realization $\mathcal{A}^{(t)}$ of the hierarchical random field, which is not contained in the training, test, or validation datasets. Moreover, we compute a Monte Carlo approximation of the expectation of the solution $\mathbb{E}[u]$ based on $100$ samples, where the individual diffusion problems are solved by FEM, PG-LOD, and our NN-LOD approach, respectively.

\subsection{Experiment 1-H: Low-contrast realizations}
\begin{figure}
     \centering
     \begin{subfigure}[a]{0.45\textwidth}
         \centering
         \includegraphics[width=\textwidth]{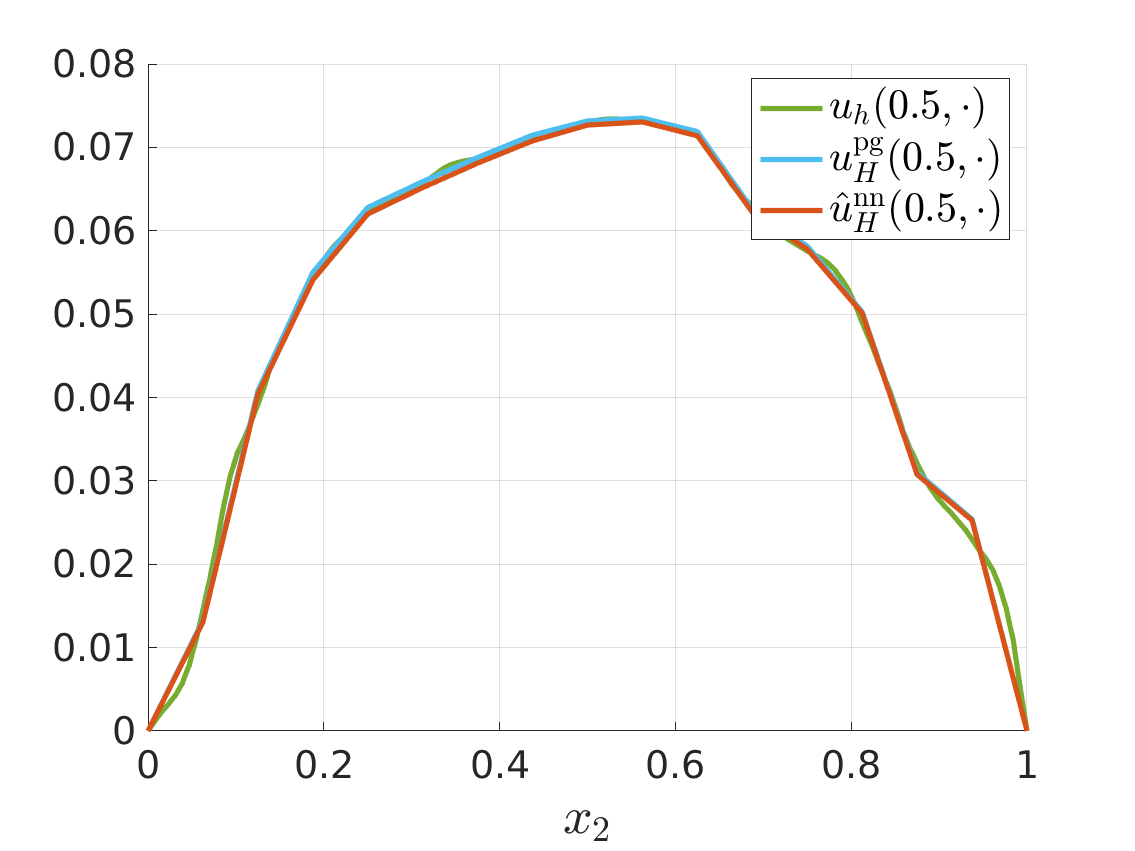}
     \end{subfigure}
     \begin{subfigure}[a]{0.45\textwidth}
         \centering
         \includegraphics[width=\textwidth]{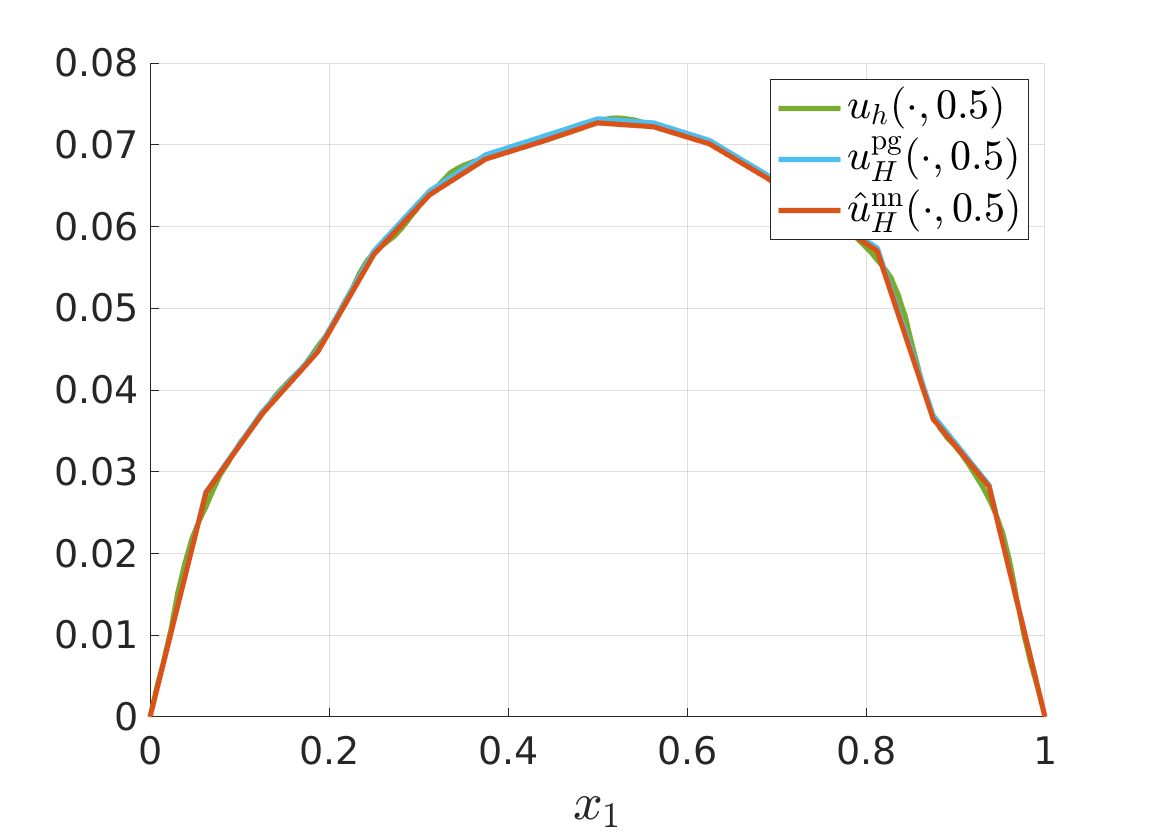}
     \end{subfigure}  
     \begin{subfigure}[a]{0.45\textwidth}
         \centering
         \includegraphics[width=\textwidth]{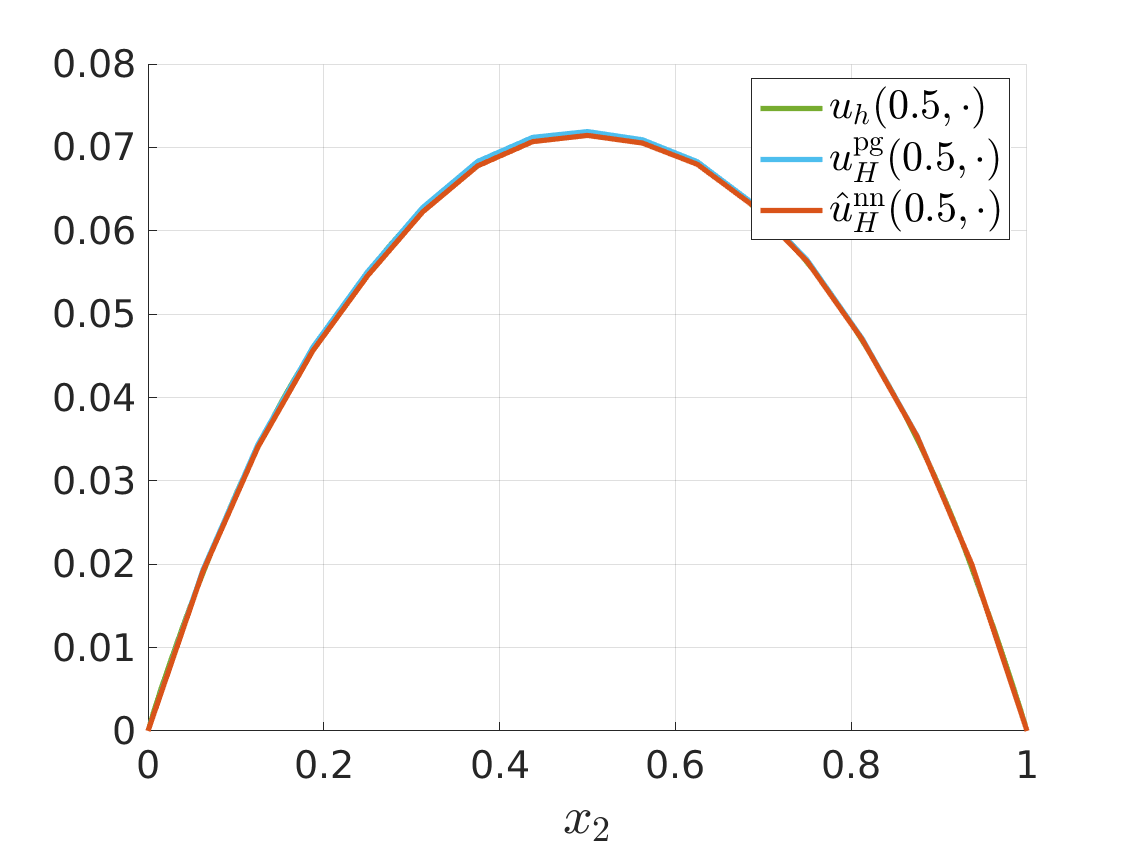}
     \end{subfigure}
     \begin{subfigure}[a]{0.45\textwidth}
         \centering
         \includegraphics[width=\textwidth]{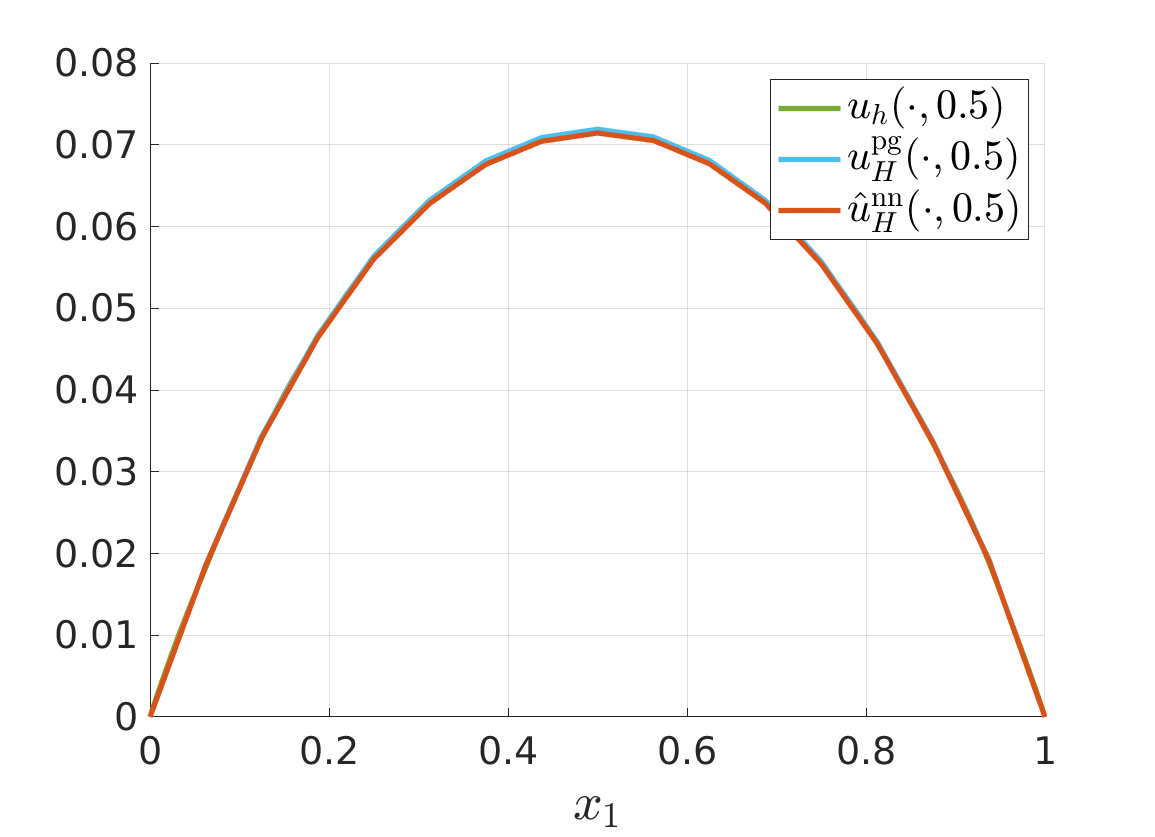}
     \end{subfigure}  
\caption{Experiment 1-H: In the top row we compare the FEM solution $u_h$ vs. the PG-LOD solution $u_H^\mathrm{pg}$ vs. the NN-LOD solution $\hat{u}_H^\mathrm{nn}$ along the cross-sections $x_1 = 0.5$ (left) and $x_2 = 0.5$ (right) for a single, fresh realization of the low-contrast hierarchical random field. 
In the bottom row we compare the Monte Carlo approximation of the mean function $\mathbb{E}[u]$ calculated with FEM samples $u_h$ vs. PG-LOD samples $u_H^\mathrm{pg}$ vs. NN-LOD samples~$\hat{u}_H^\mathrm{nn}$ along the cross-sections $x_1 = 0.5$ (left) and $x_2 = 0.5$ (right).}
\label{fig:hierarchical_loco}
\end{figure}

After training the network for $60$ epochs, we obtain a final training error of $1.42 \cdot 10^{-4}$ and a final validation error of $2.99 \cdot 10^{-4}$.
Based on the previous experiments in Section~\ref{sec:numexp:log}, we expect a good NN-LOD surrogate approximation.
Let's look at the results: The top row of Figure~\ref{fig:hierarchical_loco} compares a fine-scale FEM reference solution with the PG-LOD and the NN-based approximation for a single realization of the hierarchical random field. The bottom row of Figure~\ref{fig:hierarchical_loco} shows a Monte Carlo approximation of the expected solution based on $100$ samples for the FEM, PG-LOD, and NN-LOD discretizations, respectively.
Figure~\ref{fig:hierarchical_loco} tells us that all three spatial discretization methods (FEM, PG-LOD, NN-LOD) return qualitatively similar realizations of the solution, which are visually close in most parts of the computational domain. 
The Monte Carlo approximation of the expected solution is well resolved by both PG-LOD and NN-LOD. 

\subsection{Experiment 2-H: Moderate-contrast realizations}
\begin{figure}
     \centering
     \begin{subfigure}[a]{0.45\textwidth}
         \centering
         \includegraphics[width=\textwidth]{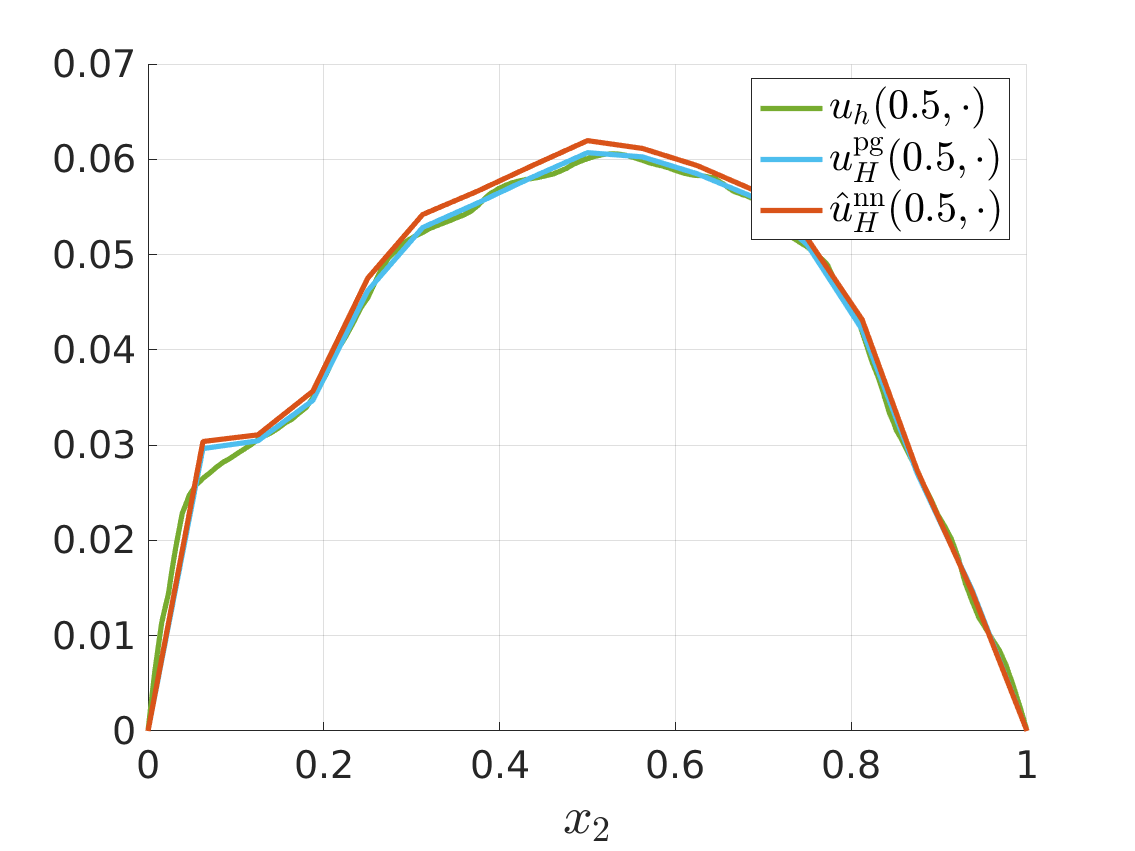}
     \end{subfigure}
     \begin{subfigure}[a]{0.45\textwidth}
         \centering
         \includegraphics[width=\textwidth]{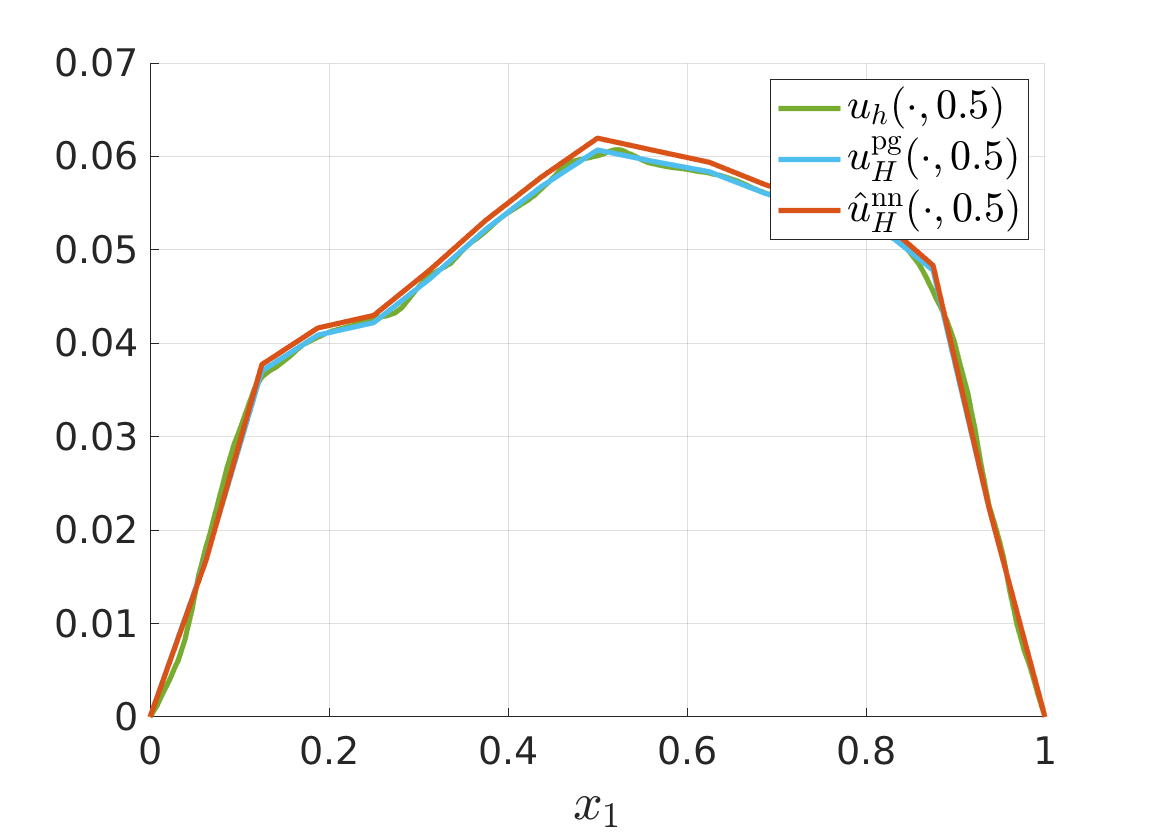}
     \end{subfigure}  
     \begin{subfigure}[a]{0.45\textwidth}
         \centering
         \includegraphics[width=\textwidth]{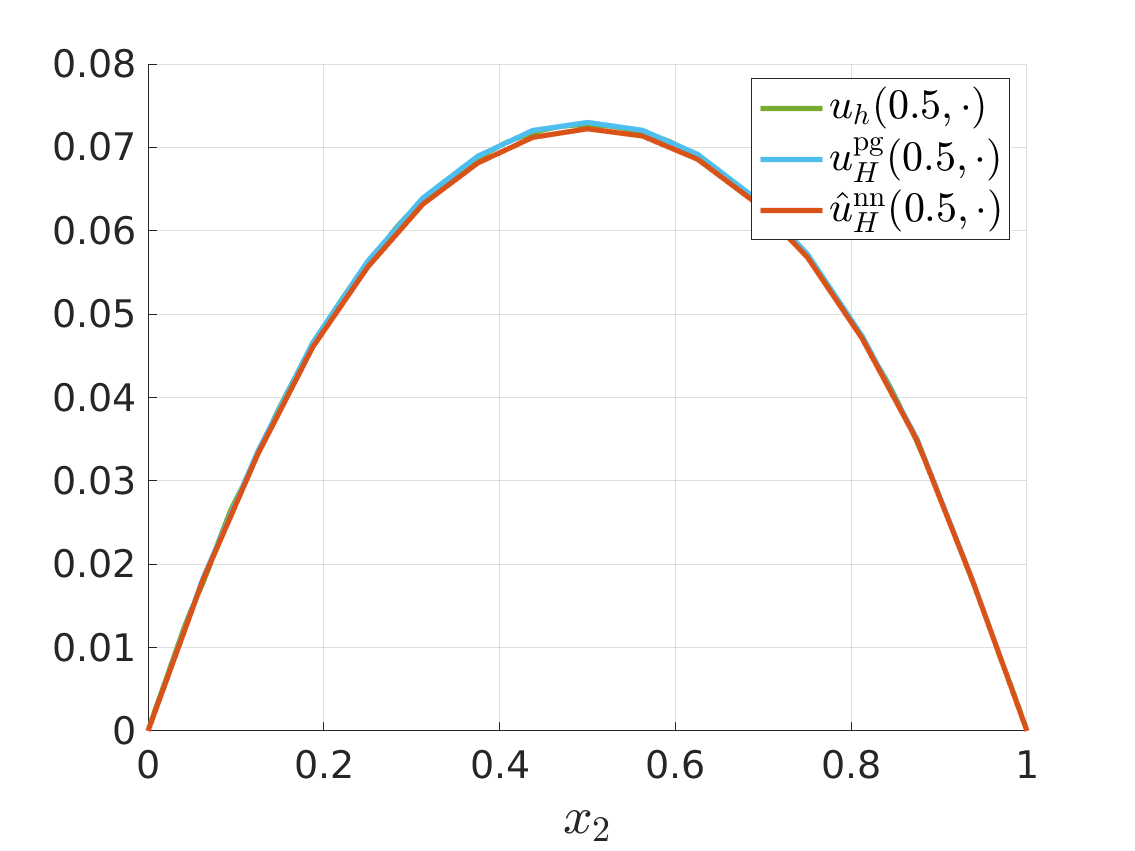}
     \end{subfigure}
     \begin{subfigure}[a]{0.45\textwidth}
         \centering
         \includegraphics[width=\textwidth]{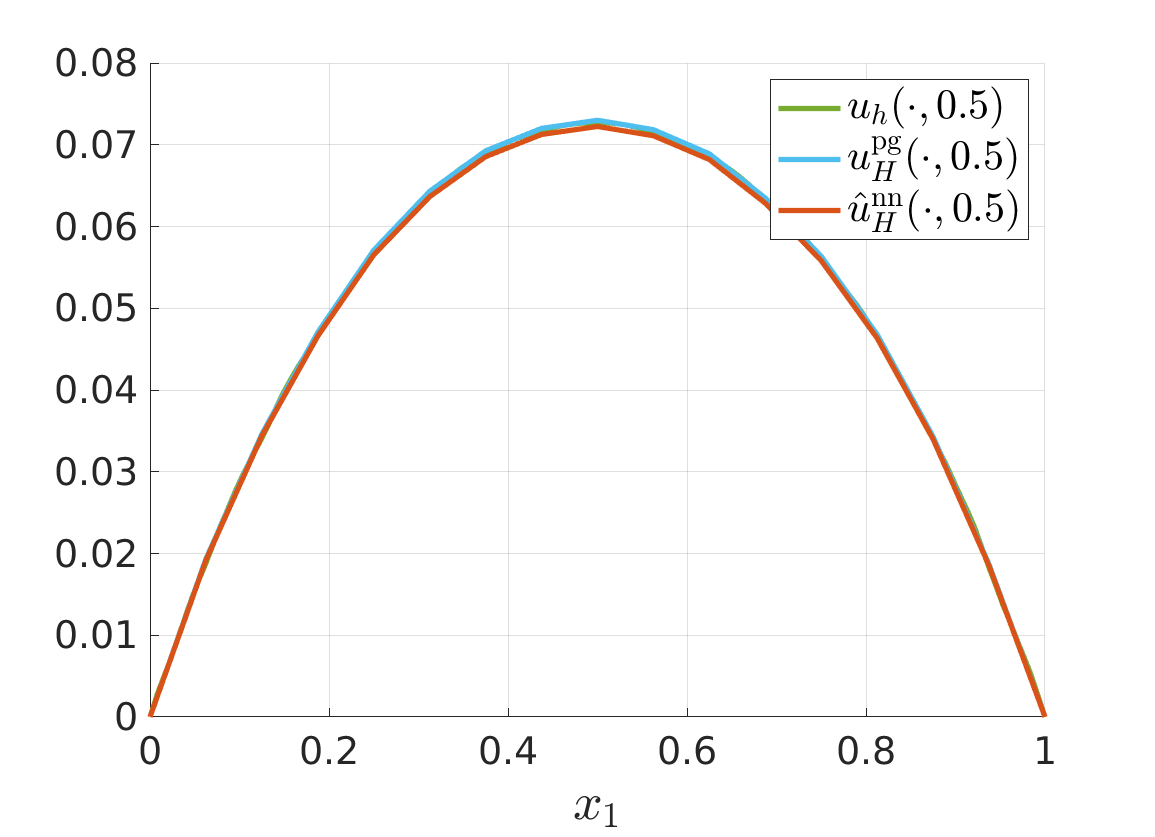}
     \end{subfigure}  
\caption{Experiment 2-H: In the top row we compare the FEM solution $u_h$ vs. the PG-LOD solution $u_H^\mathrm{pg}$ vs. the NN-LOD solution $\hat{u}_H^\mathrm{nn}$ along the cross-sections $x_1 = 0.5$ (left) and $x_2 = 0.5$ (right) for a single, fresh realization of the low-contrast hierarchical random field. 
In the bottom row we compare the Monte Carlo approximation of the mean function $\mathbb{E}[u]$ calculated with FEM samples $u_h$ vs. PG-LOD samples $u_H^\mathrm{pg}$ vs. NN-LOD samples~$\hat{u}_H^\mathrm{nn}$ along the cross-sections $x_1 = 0.5$ (left) and $x_2 = 0.5$ (right).}
\label{fig:hierarchical_mico}
\end{figure}

After training the network for $60$ epochs, we obtain a final training error of $4.25 \cdot 10^{-4}$ and a final validation error of $8.32 \cdot 10^{-4}$.
As in the low-contrast regime, we expect a good NN-LOD approximation.
A visual comparison in Figure~\ref{fig:hierarchical_mico} confirms this.
In the top row of Figure~\ref{fig:hierarchical_mico}, we see the fine-scale FEM reference solution, the PG-LOD solution, and the NN-LOD solution for a fresh realization of the hierarchical diffusion coefficient. The bottom row of Figure~\ref{fig:hierarchical_mico} shows Monte Carlo approximations of the expected solution based on $100$ samples each, where the spatial discretization is done by FEM, PG-LOD, and NN-LOD, respectively.
The conclusions are the same as those for the low-contrast realizations.

\subsection{Experiment 3-H: High-contrast realizations}
\begin{figure}
     \centering
     \begin{subfigure}[a]{0.45\textwidth}
         \centering
         \includegraphics[width=\textwidth]{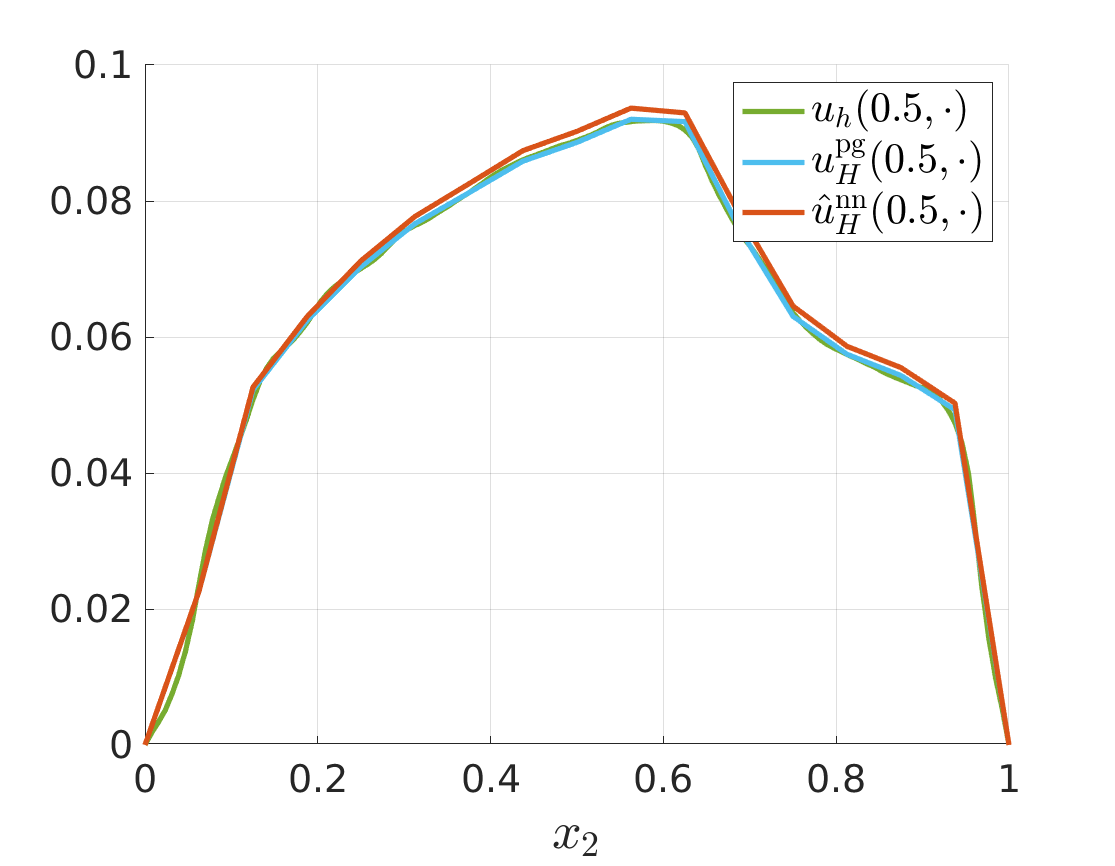}
     \end{subfigure}
     \begin{subfigure}[a]{0.45\textwidth}
         \centering
         \includegraphics[width=\textwidth]{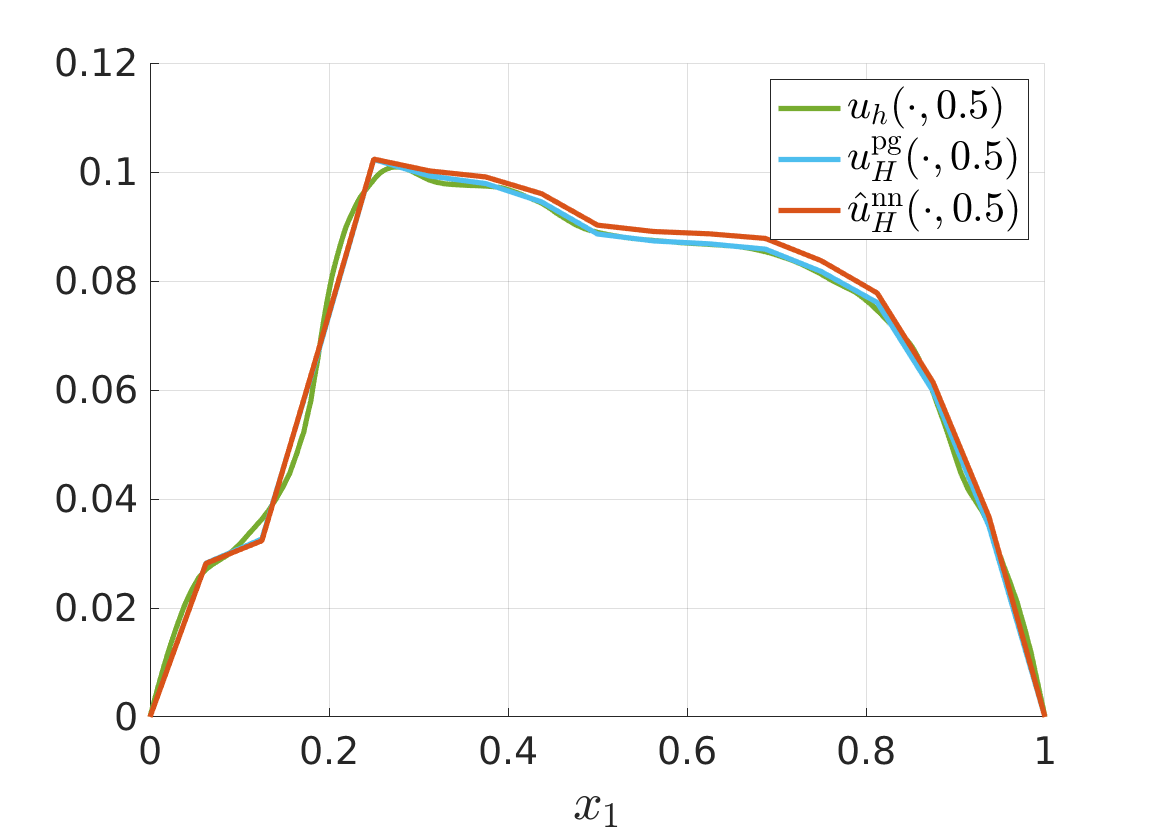}
     \end{subfigure}  
     \begin{subfigure}[a]{0.45\textwidth}
         \centering
         \includegraphics[width=\textwidth]{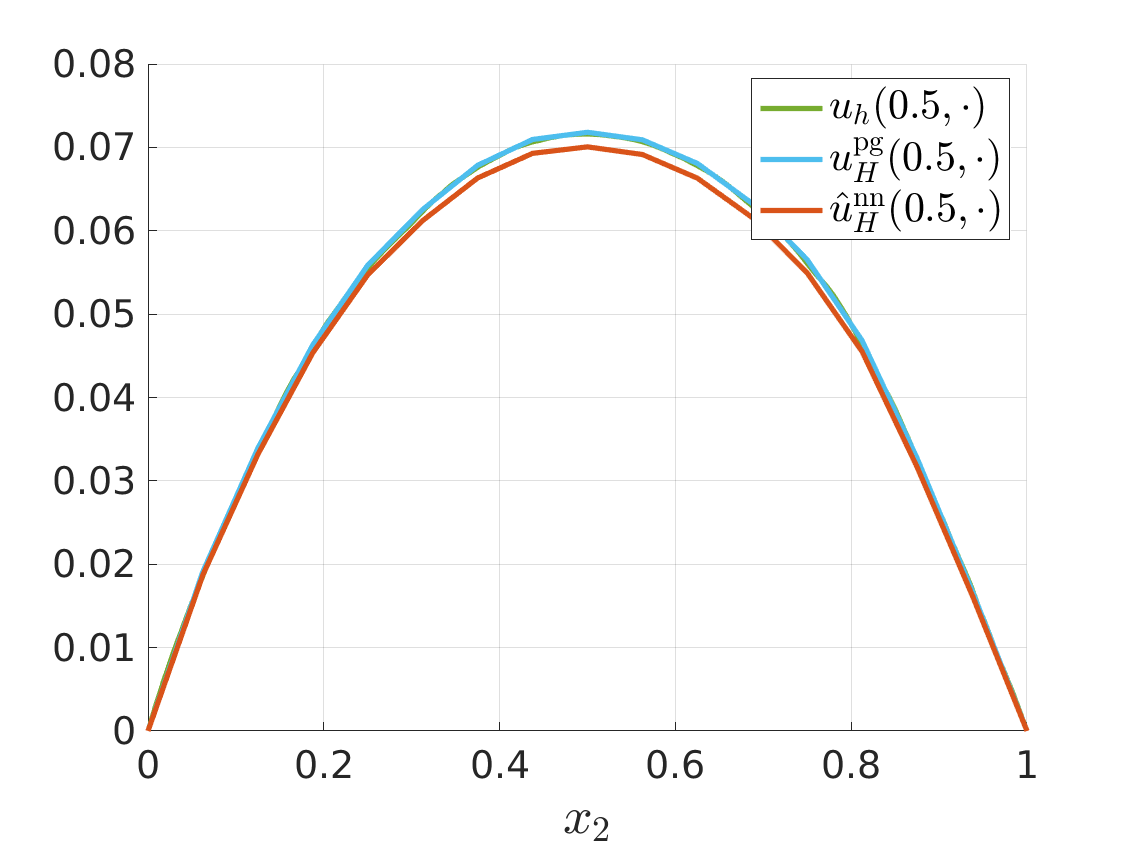}
     \end{subfigure}
     \begin{subfigure}[a]{0.45\textwidth}
         \centering
         \includegraphics[width=\textwidth]{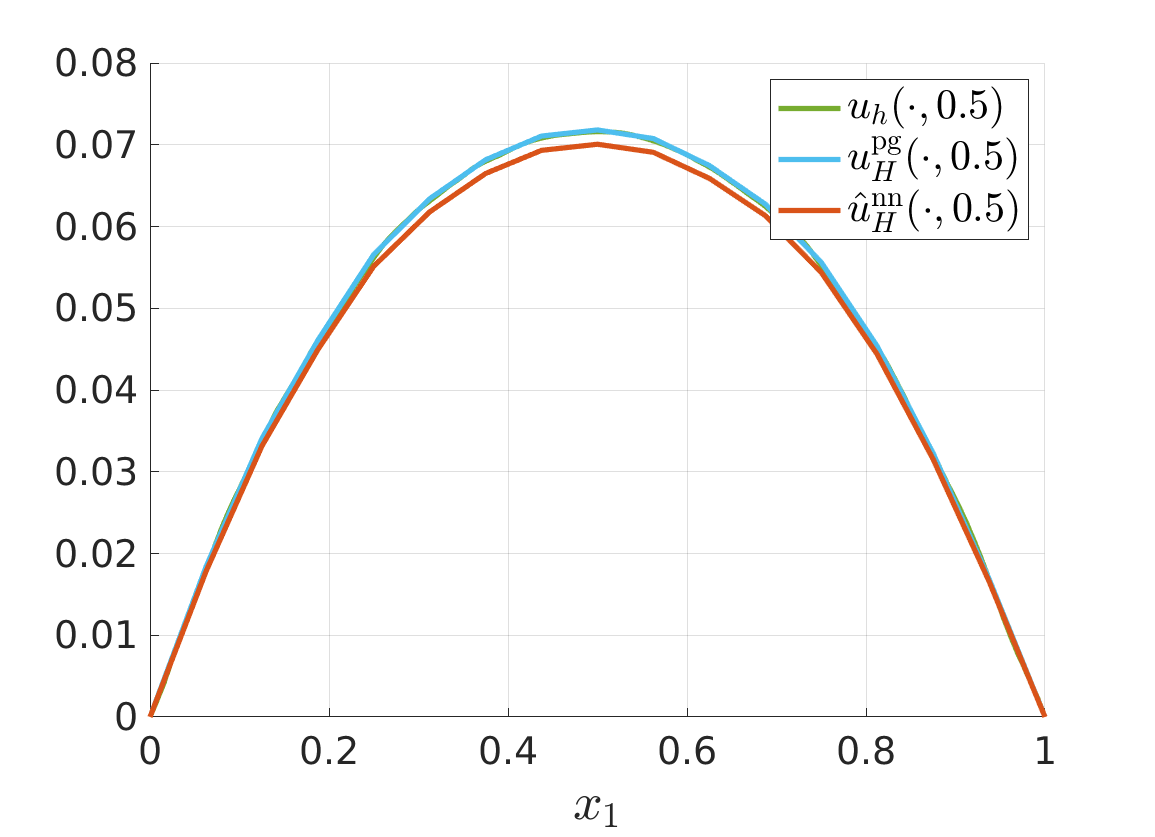}
     \end{subfigure}  
\caption{Experiment 3-H: In the top row we compare the FEM solution $u_h$ vs. the PG-LOD solution $u_H^\mathrm{pg}$ vs. the NN-LOD solution $\hat{u}_H^\mathrm{nn}$ along the cross-sections $x_1 = 0.5$ (left) and $x_2 = 0.5$ (right) for a single, fresh realization of the low-contrast hierarchical random field. 
In the bottom row we compare the Monte Carlo approximation of the mean function $\mathbb{E}[u]$ calculated with FEM samples $u_h$ vs. PG-LOD samples $u_H^\mathrm{pg}$ vs. NN-LOD samples~$\hat{u}_H^\mathrm{nn}$ along the cross-sections $x_1 = 0.5$ (left) and $x_2 = 0.5$ (right).}
\label{fig:hierarchical_hico}
\end{figure}

After training the network for $60$ epochs, we obtain a final training error of $1.9 \cdot 10^{-3}$ and a final validation error of $3.4 \cdot 10^{-3}$.
These errors are one order of magnitude larger compared to the low-contrast and moderate-contrast experiments.
Thus, we expect a decreased accuracy of the NN-LOD approximation.
Figure~\ref{fig:hierarchical_hico} compares a fine-scale FEM reference solution with the PG-LOD and the NN-based approximation for a fresh realization of the hierarchical random field. The bottom row of Figure~\ref{fig:hierarchical_hico} does the same for a Monte Carlo approximation of the expected solution based on $100$ samples each.
The individual solution approximations in the top row are qualitatively similar across all spatial discretization methods.
The Monte Carlo approximation of the expected solution based on NN-LOD samples consistently undershoots the FEM-based reference solution across the computational domain.
This is consistent with the high-contrast numerical experiments in Section~\ref{sec:numexp:log}, where we observed a similar pattern for insufficiently large datasets. 
We conjecture that this discrepancy can be reduced by enlarging the size of the training dataset for each representative correlation length.

\section{Conclusion and Open Problems}\label{sec:conclusion}
We presented a neural network--enhanced surrogate modeling approach called NN-LOD for diffusion problems with spatially varying random field coefficients with small correlation lengths. 
The method compresses fine-scale coefficients into coarse-scale surrogates using tailor-made neural networks.
The networks are trained through transfer learning, using weights and biases from surrogates obtained with uniformly bounded random field coefficients.
We carried out extensive numerical tests for lognormal diffusion coefficients that are not uniformly bounded across realizations and have varying contrasts.
These tests show that the transfer learning approach is feasible and produces NN-LOD approximations with good accuracy for low-contrast and moderate-contrast realizations.
However, the high-contrast case remains challenging since we require large training datasets to obtain good results.
As a proof-of-concept, we also carried out numerical experiments for a hierarchical diffusion coefficient with a random correlation length. 
By using a transfer learning approach with the pre-trained network surrogates from the lognormal diffusion coefficients, we obtained NN-LOD approximations with good accuracy for the low-contrast and medium-contrast realizations. 
However, in the high-contrast setting, the NN-LOD approximation shows reduced quality, which is likely caused by the limited size of the training dataset.
The precise choice of the training dataset relative to the contrast of the random field realizations is an area of future research.
Other avenues for further investigation include alternative classes of random field coefficients; for example, random fields based on Gaussian mixture models.

\bibliographystyle{plain}
\bibliography{refs}
\end{document}